\documentclass[smallextended,envcountsect]{svjour3}
\smartqed
\usepackage{lineno}
\modulolinenumbers[5]

\journalname{Journal of XXXX}
\usepackage{amssymb,amsmath,amsfonts}
\usepackage{epsfig}
\usepackage[bookmarksnumbered, colorlinks, urlcolor=blue, linkcolor=blue, citecolor=red, plainpages=false, pdfstartview=FitW]{hyperref}
\usepackage{mathptmx}
\usepackage{diagbox}
\usepackage{graphicx}
\usepackage{color}
\usepackage{caption}
\usepackage{wrapfig}
\usepackage{graphics}
\usepackage{enumerate}
\usepackage{amsxtra}
\usepackage{latexsym}
\usepackage{multirow}
\usepackage{subfigure}
\usepackage{epstopdf}
\usepackage{pdfpages}
\usepackage{algorithm}
\usepackage{algorithmic}
\usepackage{verbatim}
\usepackage{cite}
\usepackage{url}
\newcommand{\thmlist}{
\begin{list}{Step 1}
{\setlength{\leftmargin}{0.6 in}\setlength{\labelwidth} {0.5 in}}}

\newcommand{\alglist}{
\begin{list}{Step 1}
{\setlength{\leftmargin}{1.1 in} \setlength{\labelwidth}{1.0 in}}}

 \renewcommand{\proof} {\noindent {\bf Proof.} \quad}
 \newcommand{\eproof} {$\quad \square$}

\renewcommand{\subtitle}[1]{\color{blue}}

\begin{document}

\title{Residual regularization path-following methods for linear complementarity problems}
\titlerunning{Residual regularization path-following methods for LCP}
\author{Xin-long Luo \textsuperscript{$\ast$} \and Sen Zhang \and Hang Xiao}
\authorrunning{Luo \and Zhang \and Xiao}

\institute{
   Xin-long Luo, Corresponding author
   \at
   School of Artificial Intelligence,
   Beijing University of Posts and Telecommunications, P. O. Box 101,
   Xitucheng Road  No. 10, Haidian District, 100876, Beijing China,
   \email{luoxinlong@bupt.edu.cn}
   \and
   Sen Zhang
     \at
     School of Artificial Intelligence, \\
     Beijing University of Posts and Telecommunications, P. O. Box 101, \\
     Xitucheng Road  No. 10, Haidian District, 100876, Beijing China \\
     \email{senzhang@bupt.edu.cn}
     \and
   Hang Xiao
     \at
     School of Artificial Intelligence, \\
     Beijing University of Posts and Telecommunications, P. O. Box 101, \\
     Xitucheng Road  No. 10, Haidian District, 100876, Beijing China \\
     \email{xiaohang0210@bupt.edu.cn}
}

\date{Received: date / Accepted: date}
\maketitle

\begin{abstract}
In this article, we consider the residual regularization path-following method with
the trust-region updating strategy for the linear complementarity problem. This 
time-stepping selection based on the trust-region updating strategy overcomes 
the shortcoming of the line search method, which consumes the unnecessary trial 
steps in the transient-state phase. In order to improve the robustness of the 
path-following method, we use the residual regularization parameter to replace 
the traditional complementarity regularization parameter. Moreover, we prove the 
global convergence of the new method under the standard assumptions without the 
traditional assumption condition of the priority to feasibility over complementarity. 
Numerical results show that the new method is robust and efficient for the linear 
complementarity problem, especially for the dense cases. And it is more robust 
and faster than some state-of-the-art solvers such as the built-in subroutines 
PATH and MILES of the GAMS v28.2 (2019) environment. The computational time of 
the new method is about 1/3 to 1/10 of that of PATH for the dense linear 
complementarity problem.
\end{abstract}

\keywords{Continuation Newton method \and residual regularization
\and trust-region updating strategy
\and complementarity \and path-following method}

 \vskip 2mm

\subclass{90C33 \and 65K05 \and 65L05 \and 65L20}

% \linenumbers

% main text
\section{Introduction}

\vskip 2mm

In this article, we are mainly concerned with the linear complementarity problem
as follows:
\begin{align}
  y = Mx + q, \;  x_{i}y_{i} = 0, \; i = 1, \, 2, \, \ldots, \, n, \;
   x \ge 0,  \; y \ge 0,    \label{LCP}
\end{align}
where $q \in \Re^{n}$  is a vector and $M$ is an $n \times n$ positive semi-definite
matrix. For the linear complementarity problem \eqref{LCP}, there are many practical
applications such as the equilibrium of forces \cite{Erleben2007} and the economic
equilibrium problem \cite{Mathiesen1985,Rutherford1995}. And the solutions of many
problems such as the linear programming and the convex quadratic programming can be
obtained by solving it \cite{CPS2009,Wright1997}. Furthermore, there are many
efficient methods to solve it such as the Lemke's complementary pivoting algorithm
(MILES) \cite{CPS2009,LH1964,Rutherford1995}, the path-following methods
\cite{Fisher1992,Wright1994,Wright1997,Zhang1994} and their mixture method
(PATH) \cite{Dirkse1994,DF1995}.

\vskip 2mm

In this paper, we consider another path-following method based on the Newton
flow with nonnegative constraints, which is the variant of the primal-dual 
path-following method. In order to improve its robustness and efficiency, we 
use the residual regularization technique to avoid the singularity of the 
Jacobian matrix, and adopt the trust-region updating strategy to adjust the time 
step adaptively. Firstly, we construct the regularization Newton flow with nonnegative 
constraints for the linear complementarity problem \eqref{LCP} based on the 
primal-dual path-following method. Then, we use the implicit Euler method and 
the linear approximation of the quadratic function to obtain the regularization 
path-following method for following the the trajectory of the Newton flow. Finally, 
we adopt the trust-region updating strategy to adjust the time step adaptively. 
The advantage of this time-stepping selection compared with the line search method 
is that it overcomes the shortcoming of the line search method, which consumes the 
unnecessary trial steps in the transient-state phase. The other advantage is that 
it improves the robustness of the path-following method and it does not require the 
traditional assumption condition of the priority to feasibility over complementarity 
of the path-following method \cite{Wright1994,Zhang1994} when we prove its global 
convergence.

\vskip 2mm

The rest of this article is organized as follows. In the next section, we
consider the regularization path-following method and the adaptive trust-region
updating strategy for the linear complementarity problem. In section 3, we prove
the global convergence of the new method under the standard assumptions without
the traditional assumption condition of the priority to feasibility over complementarity. 
In section 4, we compare the new method with two state-of-the-art solvers, i.e.
PATH \cite{Dirkse1994,DF1995,PATH} and MILES (a Mixed Inequality and
nonLinear Equation Solver) \cite{Mathiesen1985,Rutherford1995,Rutherford2022} for
sparse problems and dense problems, test matrices of which come from the
linear programming subset of NETLIB \cite{NETLIB}. The new method is coded
with the MATLAB language and executed in MATLAB (R2020a) environment \cite{MATLAB}.
PATH and MILES are executed in the GAMS v28.2 (2019) environment \cite{GAMS}.
Numerical results show that the new method is robust and efficient for solving
the linear complementarity problem. It is more robust and faster than PATH and MILES
for the dense problems. Finally, some discussions are given in section 5.
$\|\cdot\|$ denotes the Euclidean vector norm or its induced matrix norm throughout
the paper.

\vskip 2mm

\section{Regulation path-following methods}

\subsection{The continuous Newton flow} \label{SUBSECNF}

For convenience, we rewrite the linear complementarity problem \eqref{LCP}
as the following nonlinear system of equations with nonnegativity constraints:
\begin{align}
   F(z) = \begin{bmatrix}
        y - (Mx + q)   \\
        XY e
  \end{bmatrix} = 0, \; (x, \, y) \ge 0,
  \; z = (x, \, y),  \label{NLELCP}
\end{align}
where $X = diag(x), \;  Y = diag(y)$ and all components of vector $e$ equal one.
It is not difficult to know that the Jacobian matrix $J(z)$ of $F(z)$
has the following form:
\begin{align}
  J(z) = \begin{bmatrix}
    - M & I \\
    Y &  X
  \end{bmatrix}.   \label{JZMATLCP}
\end{align}
From the second block $XYe = 0$ of equation \eqref{NLELCP}, we know that
$x_{i} = 0$ or $y_{i} = 0 \, (i = 1:n) $.  Thus, the Jacobian matrix $J(z)$
of equation \eqref{JZMATLCP} may be singular, which leads to numerical
difficulties near the solution of the nonlinear system \eqref{NLELCP}
for the Newton method or its variants. In order to overcome this difficulty,
we consider its perturbed system \cite{AG2003,Tanabe1988} as follows:
\begin{align}
  F_{\mu}(z) = F(z) - \begin{bmatrix}
                               0 \\
                               \mu e
  \end{bmatrix} = 0,
  \; z = (x, \, y) > 0, \; \mu > 0. \label{PNLEX}
\end{align}

\vskip 2mm

It is not difficult to verify that the Jacobian matrix $J(z)$ defined by
equation \eqref{JZMATLCP} is nonsingular when $M$ is a positive semi-definite
matrix and $(x, \; y) > 0$ (see Lemma 5.9.8, p. 469 in \cite{CPS2009}). Thus, by
using the implicit theorem and the inequality $(a - b)(c-d) \le |ac - bd|
\; (a > 0, \, b > 0, \, c > 0, \, d > 0)$, the perturbed system \eqref{PNLEX}
has a unique solution when $M$ is a positive definite matrix (see Theorem 5.9.13,
p. 471 in \cite{CPS2009}) for its detailed proof). The solution
$z(\mu)$ of the perturbed system \eqref{PNLEX} defines the central path, and
$z(\mu)$ approximates the solution $z^{\ast}$ of the nonlinear system
\eqref{NLELCP} when $\mu$ tends to zero \cite{Wright1997,CPS2009}.

\vskip 2mm

If the damped Newton method is applied to the perturbation system \eqref{PNLEX}
\cite{Kelley2003,NW1999}, we have
\begin{align}
     z_{k+1} = z_{k} - \alpha_{k} J(z_{k})^{-1}  F_{\mu}(z_{k}), \label{NEWTON}
\end{align}
where $J(z_{k})$ is the Jacobian matrix of $F_{\mu}(z)$. We regard
$z_{k+1} = z(t_{k} + \alpha_{k}), \; z_{k} = z(t_{k})$ and let
$\alpha_{k} \to 0$, then we obtain the continuous Newton flow with nonnegativity
constraints \cite{AS2015,Davidenko1953,LY2021,PHP1975,Tanabe1988} of the
perturbed system \eqref{PNLEX} as follows:
\begin{align}
  \frac{dz(t)}{dt} = - J(z)^{-1}F_{\mu}(z), \hskip 2mm  z = (x, \, y) > 0.
        \label{NEWTONFLOW}
\end{align}
Actually, if we apply an iteration with the explicit Euler method
\cite{SGT2003} for the continuous Newton flow \eqref{NEWTONFLOW}, we
obtain the damped Newton method \eqref{NEWTON}.

\vskip 2mm

Since the Jacobian matrix $J(z) = F'_{\mu}(z)$ may be singular, we reformulate
the continuous Newton flow \eqref{NEWTONFLOW} as the following general formula
for the linear complementarity problem \eqref{NLELCP} :
\begin{align}
    & -M \frac{dx(t)}{dt} + \frac{dy(t)}{dt} = - r_{q}(x, \, y), \label{LFRQDAE} \\
    & Y \frac{dx(t)}{dt} + X \frac{dy(t)}{dt} = - (XYe - \mu(t) e), \;  (x, \, y) > 0,
    \label{DAEFLOW}
\end{align}
where the residual $r_{q}(x, \, y) = y - (Mx+q)$.
The continuous Newton flow \eqref{LFRQDAE}-\eqref{DAEFLOW} has some nice
properties. We state one of them as the following property \ref{PRODAEFLOW}
\cite{Branin1972,LXL2022,LY2021,LX2021,Tanabe1979}.

\vskip 2mm

\begin{property} \label{PRODAEFLOW}
Assume that  $(x(t), \, y(t))$ is the solution of the
continuous Newton flow \eqref{LFRQDAE}-\eqref{DAEFLOW}, then
$r_{q}(x(t), \, y(t))$ converges to zero and
$x_{i}(t)y_{i}(t) \, (i = 1, \, 2, \, \ldots, \, n)$ converge to zero when
$0 \le \mu(t) \le \sigma \min_{1\le i \le n}\{x_{i}(t)y_{i}(t)\} \; (0 < \sigma < 1)$
and $t \to \infty$. That is to say, for every limit point $(x^{\ast}, \, y^{\ast})$
of $(x(t), \, y(t))$, it is also a solution of the linear complementarity problem
\eqref{NLELCP}. Furthermore, $x(t)$ and $y(t)$ keep positive values when the initial
point $\left(x^{0}_{i}, \; y^{0}_{i}\right) > 0 \, (i = 1, \, 2, \, \ldots, \, n)$.
\end{property}
\proof Assume that $z(t)$ is the continuous solution of the continuous Newton flow
\eqref{LFRQDAE}-\eqref{DAEFLOW}, then we have
\begin{align}
    & \frac{d}{dt} r_{q}(x, \, y) = -M \frac{dx}{dt} + \frac{dy}{dt}
     = - r_{q}(x,\, y), \label{ODELEQ} \\
    & \frac{d}{dt} (XYe) = X \frac{dy}{dt} + Y \frac{dx}{dt}
    = - (XYe - \mu(t) e). \label{ODECEQ}
\end{align}
Consequently, from equations \eqref{ODELEQ}-\eqref{ODECEQ} and
$0 \le \mu(t) \le \sigma \min_{1\le i \le n}\{x_{i}(t)y_{i}(t)\}$, we obtain
\begin{align}
     & r_{q}(x(t), \, y(t)) = r_{q}\left(x^{0}, \, y^{0}\right) \exp(-t),
      \label{FUNPARQ} \\
     &  -XYe \le \frac{d}{dt} (XYe) \le -(1-\sigma)XYe. \label{FUNPAR}
\end{align}
From equations \eqref{FUNPARQ}-\eqref{FUNPAR}, we know that
$r_{q}(x(t), \, y(t))$ converges to zero with the linear rate of convergence
when $t$ tends to infinity. Furthermore, from equation \eqref{FUNPAR} and
the Gronwall inequality \cite{Gronwall1919}, we have
\begin{align}
  x_{i}^{0}y_{i}^{0} \exp(-t) \le x_{i}(t)y_{i}(t)
  \le x_{i}^{0}y_{i}^{0} \exp(-(1 - \sigma)t), \; i = 1, \, 2, \, \ldots, \, n.
  \label{SCGINEQ}
\end{align}
Consequently, from equation \eqref{SCGINEQ}, we know that
$x_{i}(t)y_{i}(t) \ge 0 \; (i = 1, \, 2, \, \ldots, \, n)$ and
$\lim_{t \to \infty} x_{i}(t)y_{i}(t) = 0 \; (i = 1, \, 2, \, \ldots, \, n)$.

\vskip 2mm

From equation \eqref{SCGINEQ}, we know that $x_{i}(t)y_{i}(t) > 0 \; (t \ge 0)$
when $\left(x^{0}, \; y^{0}\right) > 0$. If we have $x_{i}(T) < 0$ or
$y_{i}(T) < 0$ for a finite value $T > 0$, there exists $\bar{t} \in (0, \; T)$
such that $x_{i}(\bar{t}) = 0$ or $y_{i}(\bar{t}) = 0$. Consequently, we have
$x_{i}(\bar{t})y_{i}(\bar{t}) = 0$, which contradicts $x_{i}(\bar{t})y_{i}(\bar{t}) > 0$.
Thus, we have $(x(t), \, y(t)) > 0 $ for all $t > 0$. Therefore, if the
solution $(x(t), \, y(t))$ of the continuous Newton flow
\eqref{LFRQDAE}-\eqref{DAEFLOW} belongs to a compact set, there exists a limit point
$(x^{\ast}, \, y^{\ast})$ when $t$ tends to infinity, and this limit point
$(x^{\ast}, \, y^{\ast})$ is also a solution of the linear complementarity problem
\eqref{NLELCP}.  \qed

\vskip 2mm

\begin{remark} The inverse $J(z)^{-1}$ of the Jacobian matrix  $J(z)$ can be
regarded as the pre-conditioner of $F_{\mu}(z)$  such that the every element
$z_{i}(t)$ of $z(t)$ has roughly the same rate of convergence and it mitigates
the stiffness of the ODE \eqref{NEWTONFLOW} \cite{LXL2022,LX2021}. This property is
very useful since it makes us adopt the explicit ODE method to follow the trajectory
of the Newton flow \eqref{NEWTONFLOW} efficiently.
\end{remark}

\vskip 2mm

\subsection{The residual regularization path-following method} \label{SUBSECRPFM}

\vskip 2mm

From property \ref{PRODAEFLOW}, we know that the continuous Newton flow
\eqref{LFRQDAE}-\eqref{DAEFLOW} has the global convergence. However, when
the Jacobian matrix $J(z)$ is singular or nearly singular, the ODE
\eqref{LFRQDAE}-\eqref{DAEFLOW} is the system of differential-algebraic equations
\cite{AP1998,BCP1996,HW1996} and its trajectory can not be efficiently solved
by the general ODE method \cite{BJ1998,JT1995} such as the backward differentiation
formulas (the built-in subroutine ode15s.m of the MATLAB environment \cite{MATLAB,SGT2003}).
Thus, we need to construct the special method to solve this problem. Furthermore,
we expect that the new method has the global convergence as the homotopy
continuation methods \cite{AG2003,OR2000} and the fast rate of convergence as the
traditional optimization methods. In order to achieve these two aims, we consider
the continuation Newton method and the trust-region updating strategy for problem
\eqref{LFRQDAE}-\eqref{DAEFLOW}.

\vskip 2mm

We apply the implicit Euler method \cite{AP1998,BCP1996} to the continuous Newton
flow \eqref{LFRQDAE}-\eqref{DAEFLOW}, then we obtain
\begin{align}
    & -M \frac{x^{k+1} - x^{k}}{\Delta t_{k}} + \frac{y^{k+1} - y^{k}}{\Delta t_{k}}
    = - r_{q}\left(x^{k+1}, \, y^{k+1}\right), \label{IMLEQ} \\
    & Y^{k+1} \frac{x^{k+1} - x^{k}}{\Delta t_{k}}
    + X^{k+1} \frac{y^{k+1} - y^{k}}{\Delta t_{k}}
    = - \left(X^{k+1}Y^{k+1}e - \mu(t_{k+1}) e\right).   \label{IMEDAE}
\end{align}
Since equation \eqref{IMEDAE} is a nonlinear system, it is not directly solved,
and we seek for its explicit approximation formula. We replace $Y^{k+1}$
and $X^{k+1}$ with $Y^{k}$ and $X^{k}$ in the left-hand side of equation
\eqref{IMEDAE}, respectively. And we substitute $X^{k+1}Y^{k+1}e$ with its
linear approximation $X^{k}Y^{k}e + \frac{\Delta t_{k}}{1 + \Delta t_{k}}
(Y^{k}\Delta x^{k} + X^{k}\Delta y^{k})$ in the right-hand side of equation
\eqref{IMEDAE}. We set $\mu(t_{k+1}) = \sigma_{k} \mu_{k}$. Then, we obtain the
continuation Newton method (one of path-following methods) as follows:
\begin{align}
     & -M \Delta x^{k} + \Delta y^{k} = - r_{q}^{k}, \label{CONNML} \\
     & Y^{k} \Delta x^{k} + X^{k}\Delta y^{k} = - r_{c}^{k}, \label{CONNMC} \\
     & x^{k+1} = x^{k} + \frac{\Delta t_k}{1+\Delta t_k} \Delta x^{k}, \;
     y^{k+1} = y^{k} + \frac{\Delta t_k}{1+\Delta t_k} \Delta y^{k},
     \label{CONNMZ1}
\end{align}
where $r_{q}^{k} = r_{q}(x^{k}, \, y^{k})$, $r_{c}^{k} = X^{k}y^{k}
- \sigma_{k} \mu_{k} e$, $0 < \sigma_{min} \le \sigma_{k} \le \sigma_{max} < 1$
and the perturbation parameter $\mu_{k}$ is selected as follows:
\begin{align}
     \mu_{k} = \frac{(x^{k})^{T}y^{k} + \|r_{q}^{k}\|}{2n}. \label{DEFMUK}
\end{align}
The selection of $\mu_{k}$ in equation \eqref{DEFMUK} is slightly different
to the traditional selection $\mu_{k} = \left(x^{k}\right)^{T}y^{k}/n$
\cite{Wright1994,Wright1997,Zhang1994}. According to our numerical
experiments, the selection of $\mu_{k}$ in equation \eqref{DEFMUK}
can improve the robustness of the path-following method, in comparison to
the traditional selection selection $\mu_{k} = \left(x^{k}\right)^{T}y^{k}/n$.

\vskip 2mm

\begin{remark} If we set $\alpha_{k} = \Delta t_k/(1+\Delta t_k)$ in equation
\eqref{CONNMZ1}, we obtain the damped Newton method (the primal-dual path-following
method) \eqref{NEWTON}. However, from the view of the ODE method, they are different.
The damped Newton method \eqref{NEWTON} is derived from the explicit Euler method
applied to the continuous Newton flow \eqref{LFRQDAE}-\eqref{DAEFLOW}. Its time
step $\alpha_k$ is restricted by the numerical stability \cite{HW1996,SGT2003}.
That is to say, for the linear test equation $dx/dt = - \lambda x \; (\lambda > 0)$,
its time step $\alpha_{k}$ is restricted by the stable region
$|1-\lambda \alpha_{k}| \le 1$. Therefore, the large step $\alpha_{k}$ can not
be adopted in the steady-state phase. The continuation Newton method
\eqref{CONNML}-\eqref{CONNMZ1} is derived from the implicit Euler method applied
to the continuous Newton flow \eqref{LFRQDAE}-\eqref{DAEFLOW} and the linear
approximation of $F_{\mu(t_{k+1})}(z_{k+1})$, and its time step $\Delta t_k$ is
not restricted by the numerical stability for the linear test equation. Therefore,
the large time step $\Delta t_{k}$ can be adopted in the steady-state
phase, and the continuation Newton method \eqref{CONNML}-\eqref{CONNMZ1}
mimics the Newton method. Consequently, it has the fast rate of convergence near
the solution $z^{\ast}$ of the nonlinear system \eqref{NLELCP}. The most of all,
the substitution $\alpha_{k}$ with $\Delta t_{k}/(\Delta t_{k} + 1)$ is
favourable to adopt the trust-region updating strategy for adaptively adjusting
the time step $\Delta t_{k}$ such that the continuation Newton method
\eqref{CONNML}-\eqref{CONNMZ1} accurately follows the trajectory of the
continuation Newton flow \eqref{LFRQDAE}-\eqref{DAEFLOW} in the transient-state
phase and achieves the fast rate of convergence in the steady-state phase.
\end{remark}

\vskip 2mm

When the diagonal matrix $X^{k}$ is positive definite, from equations
\eqref{CONNML}-\eqref{CONNMC}, $\Delta x^{k}$ and $\Delta y^{k}$
can be solved by the following two subsystems:
\begin{align}
   & \left(M + (X^{k})^{-1}Y^{k}\right) \Delta x^{k}
   = r_{q}^{k} - (X^{k})^{-1}r_{c}^{k},    \label{DELTXK} \\
   & \Delta y^{k} = M \Delta x^{k} - r_{q}^{k}.   \label{DELTYK}
\end{align}
When matrix $M$ is positive semi-definite and $(x^{k}, \, y^{k}) > 0$, the
left hand-side matrix is positive definite. Thus, the system \eqref{DELTXK}
can be solved by the partial pivoting Gaussian elimination method
(see pp. 125-130, \cite{GV2013}).

\subsection{The time-stepping control and the initial point selection}

\vskip 2mm

Another issue is how to adaptively adjust the time step $\Delta t_k$
at every iteration. A popular and efficient time-stepping control is based on
the trust-region updating strategy
\cite{CGT2000,Deuflhard2004,Higham1999,Luo2010,Luo2012,LXLZ2021,LX2021,LX2021B,
LX2022,LY2021,LLS2022,LXL2022,Yuan2015}. Its main idea can be described as
follows. Firstly, we construct a merit function reflecting the feasibility
$r_{q}(x, \, y) = 0$ and the complementarity
$x_{i}y_{i} = 0 \; (i = 1, \, 2, \, \ldots, \, n)$ as
\begin{align}
    \phi(x, \, y) = x^{T}y + \|r_{q}(x, \, y)\|,   \label{MERITFUN}
\end{align}
where $r_{q}(x, \, y) = y - (Mx + q)$ and $(x, \, y) \ge 0$.

\vskip 2mm

Then, we define the linear approximation model $m_{k}$ of
$\phi(x^{k} + \alpha_{k} \Delta x^{k}, \, y^{k} + \alpha_{k} \Delta y^{k})$
as
\begin{align}
    m_{k}(\alpha) & = \phi(x^{k}, \, y^{k}) +
    \alpha \left((y^{k})^{T}\Delta x^{k} + (x^{k})^{T} \Delta y^{k}
     - \|r_{q}^{k}\|\right),     \label{LMAPP}
\end{align}
where $\alpha = \Delta t /(1+\Delta t)$ and $(\Delta x^{k}, \, \Delta y^{k})$
is the Newton step and solved by equations \eqref{CONNML}-\eqref{CONNMC}.

\vskip 2mm

Finally, we adaptively adjust the time step $\Delta t_{k}$ according to the
difference between $m_{k}(\alpha_{k})$ and $\phi(x^{k}
+ \alpha_{k} \Delta x^{k}, \, y^{k} + \alpha_{k} \Delta y^{k})$.
Namely, the time step $\Delta t_{k+1}$ will be enlarged
when $m_{k}(\alpha_{k})$ approximates $\phi(x^{k} +
\alpha_{k} \Delta x^{k}, \, y^{k} + \alpha_{k} \Delta y^{k})$
well, and $\Delta t_{k+1}$ will be reduced when $m_{k}(\alpha_{k})$
approximates $\phi(x^{k} + \alpha_{k} \Delta x^{k}, \, y^{k}
+ \alpha_{k} \Delta y^{k})$ badly.

\vskip 2mm

We define the predicted reduction
$Pred_{k}$ and the actual reduction $Ared_{k}$ as follows:
\begin{align}
    & Pred_{k} = \phi(x^{k}, \, y^{k})  - m_{k}(\alpha_{k})
    = \alpha_{k} \left(\|r_{q}^{k}\| -(y^{k})^{T}\Delta x^{k}
    - (x^{k})^{T} \Delta y^{k}\right), \label{PREDK} \\
    & Ared_{k} = \phi(x^{k}, \, y^{k})
    - \phi(x^{k} + \alpha_{k} \Delta x^{k},
    \, y^{k} + \alpha_{k} \Delta y^{k}) \nonumber \\
    & \hskip 2mm = \alpha_{k} \left(\|r_{q}^{k}\| -(y^{k})^{T}\Delta x^{k}
    - (x^{k})^{T} \Delta y^{k}\right)
    - \alpha_{k}^{2} (\Delta x^{k})^{T}\Delta y^{k}. \label{AREDK}
\end{align}
Then, we enlarge or reduce the time step $\Delta t_{k+1}$ at every iteration according
to the following ratio:
\begin{align}
  \rho_k = \frac{Ared_{k}}{Pred_{k}}
  = \frac{\|r_{q}^{k}\| -(y^{k})^{T}\Delta x^{k} - (x^{k})^{T} \Delta y^{k}
    - \alpha_{k} (\Delta x^{k})^{T}\Delta y^{k}}
    {\|r_{q}^{k}\| -(y^{k})^{T}\Delta x^{k}
    - (x^{k})^{T} \Delta y^{k}}, \label{RHOK}
\end{align}
where $\alpha_{k} = \Delta t_{k}/(1+\Delta t_{k})$. A particular adjustment
strategy is given as follows:
\begin{align}
   \Delta t_{k+1} =
     \begin{cases}
          2 \Delta t_k, \;  \text{if} \; \rho_{k} \ge \eta_{2} \;
          \text{and} \; (x^{k+1}, \, y^{k+1}) > 0, \\
         \Delta t_k, \; \text{else if} \; \eta_1 \le \rho_{k} < \eta_{2} \;
         \text{and} \; (x^{k+1}, \, y^{k+1}) > 0, \\
    \frac{1}{2}  \Delta t_k, \; \text{others},
    \end{cases} \label{TSK1}
\end{align}
where the constants $\eta_{1}, \; \eta_{2}$ are selected as
$\eta_1 = 0.25, \; \eta_2 = 0.75$, respectively. We set
\begin{align}
    (x^{k+1}, \, y^{k+1}) = (x^{k}, \, y^{k})
    + \frac{\Delta t_{k}}{1+\Delta t_{k}}
    (\Delta x^{k}, \, \Delta y^{k}).   \label{ACCPXK1}
\end{align}
When $\rho_{k} \ge \eta_{a}$ and $(x^{k+1}, \, y^{k+1}) > 0$,
we accept the trial step, otherwise we discard the trial step and set
\begin{align}
   (x^{k+1}, \, y^{k+1}) = (x^{k}, \, y^{k}), \label{NOACXK1}
\end{align}
where $\eta_{a}$ is a small positive number such as $\eta_{a} = 1.0\times 10^{-6}$.

\vskip 2mm

\begin{remark}
This new time-stepping control based on the trust-region updating strategy
has some advantages compared to the traditional line search strategy.
If we use the line search strategy and the damped Newton method \eqref{NEWTON}
to follow the trajectory $z(t)$ of the continuous Newton flow \eqref{NEWTONFLOW},
in order to achieve the fast rate of convergence in the steady-state phase,
the time step $\alpha_{k}$ of the damped Newton method is tried from 1
and reduced by the half with many times at every iteration. Since the linear
model $F_{\sigma_{k} \mu_{k}}(z_{k}) + J(z_{k})\Delta z_{k}$ may not
approximate $F_{\sigma_{k}\mu_{k}}(z_{k}+\Delta z_{k})$ well in the transient-state phase,
the time step $\alpha_{k}$ will be small. Consequently, the line search
strategy consumes the unnecessary trial steps in the transient-state phase.
However, the selection of the time step $\Delta t_{k}$ based on the
trust-region updating strategy \eqref{RHOK}-\eqref{TSK1} can overcome this shortcoming.
\end{remark}

\vskip 2mm

In order to ensure that the algorithm works well for the general linear
complementarity problem, the initial point selection is also a key point. We select
the starting point $(x^{0}, \, y^{0})$ as follows:
\begin{align}
   x^{0} = 10*e, \; v^{0} = Mx^{0} + q, \; y^{0}_{i}
    = \begin{cases}
          v_{i}^{0}, \; \text{if} \; v_{i}^{0} > 0,  \\
          10^{-3}, \; \text{otherwise}.
    \end{cases} \label{INIPTS}
\end{align}

\vskip 2mm

In order to improve the stability of the algorithm, we add a small regularization
item $\upsilon I$ to matrix $M$ \cite{CPS2009,Gana1982,Hansen1994,Venkateswaran1993}
when $\mu_{k}$ is not small, where $\mu_{k}$ is defined
by equation \eqref{DEFMUK}.
Specifically, we adopt the following strategy as the regularizer $M_{\upsilon}$ of matrix
$M$:
\begin{align}
     M_{\upsilon} = \begin{cases}
             M + \upsilon I, \; \text{if} \; \mu_{k} \ge  \upsilon, \\
             M, \; \text{otherwise},
         \end{cases} \label{REGMAT}
\end{align}
where $\upsilon = 10^{-3}$.

\vskip 2mm

According to the above discussions, we give the detailed descriptions of the
path-following method and the trust-region updating strategy for
the monotone linear complementarity problem \eqref{LCP} in Algorithm \ref{ALGRPFM}.

\vskip 2mm

\begin{algorithm}
    \renewcommand{\algorithmicrequire}{\textbf{Input:}}
	\renewcommand{\algorithmicensure}{\textbf{Output:}}
	\caption{Regularization path-following methods with the trust-region updating strategy for
    linear complementarity problems (The RPFMTr method)}
    \label{ALGRPFM}
    \begin{algorithmic}[1]
        \REQUIRE ~~\\
          matrix $M \in \Re^{n \times n}$ and vector $q \in \Re^{n}$
        for the problem:
        $y = Mx + q, \; x_{i}y_{i} = 0 \; (i = 1:n), \; x \ge 0, \; y \ge 0$.
		\ENSURE ~~ \\
        the linear complementarity solution: $(solx, \, soly)$.
        \STATE Initialize parameters: $\eta_{a} = 10^{-6}$, $\eta_1 = 0.25$,
        $\eta_2 = 0.75$, $\epsilon = 10^{-6}$, $\Delta t_0 = 10^{-2}$,
        bigMfac = 10, $\upsilon = 10^{-3}$, $\sigma_{0} = 0.5$, maxit = 600.
        \STATE Initialize $x^0$ = bigMfac*ones(n, 1);  $y^{0} = Mx^{0} + q$;
        $y^{0}(y^{0} < 0) = 10^{-3}$.
        \STATE Regularize matrix: $M_{\upsilon} = M + \upsilon I$.
        \STATE Set flag\_success\_trialstep = 1, $\text{itc} = 0, \; k = 0$.
        \WHILE {(itc $<$ maxit)}
           \IF{(flag\_success\_trialstep == 1)}
              \STATE Set itc = itc + 1.
              \STATE  Compute $r_{q}^{k} = y^{k} - (M_{\upsilon}x^{k} + q)$;
              $\mu_{k} = (\|r_{q}^{k}\| + (x^{k})^{T}y^{k})/(2n)$.
              \STATE  Compute $\text{Resk} = \max \{\|x^{k}.y^{k}\|_{\infty}, \;
              \|y^{k} - (Mx^{k} + q)\|_{\infty}\}$.
              \IF{($\text{Resk} < \epsilon$)}
                 \STATE break;
              \ENDIF
              \STATE Set $\sigma_k = \min \{\sigma_{k}, \; \mu_k\}$.
              \STATE Compute $r_{c}^{k} = x^{k}.y^{k} - \sigma_{k}\mu_{k}*\text{ones}(n,1)$.
              \STATE By solving the linear system \eqref{DELTXK}-\eqref{DELTYK},
              we obtain $\Delta x^{k}$ and $\Delta y^{k}$.
           \ENDIF
           \STATE Set $(x^{k+1}, \, y^{k+1}) = (x^{k}, \, y^{k})
              + \frac{\Delta t_{k}}{1 + \Delta t_{k}} (\Delta x^{k}, \, \Delta y^{k})$.
           \STATE Compute the ratio $\rho_{k}$ from equation \eqref{RHOK} and adjust
           $\Delta t_{k+1}$ according to the formula \eqref{TSK1}.
           \IF{(($\rho_{k} \ge \eta_{a}$) \&\& ($x^{k+1}, \, y^{k+1}) > 0$))}
               \STATE Accept the trial point $(x^{k+1}, \, y^{k+1})$;
               Set flag\_success\_trialstep = 1.
               \IF{($\|x^{k+1} - x^{k}\|_{\infty} > 0.1$)}
                  \STATE Set $\sigma_{k+1} = 0.5$.
               \ELSE
                  \STATE Set $\sigma_{k+1} = 0.1$.
               \ENDIF
               \IF{($\mu_{k} < \upsilon$)}
                   \STATE Set $M_{\upsilon} = M$.
               \ENDIF
           \ELSE
               \STATE Set $(x^{k+1}, \, y^{k+1}) = (x^{k}, \, y^{k})$;
               flag\_success\_trialstep = 0.
           \ENDIF
           \STATE Set $k \leftarrow k+1$.
        \ENDWHILE
        \STATE Return $(solx, \, soly) =  (x^{k}, \, y^{k})$.
   \end{algorithmic}
\end{algorithm}

\vskip 2mm

\section{Convergence analysis}

\vskip 2mm

In this section, we analyze the global convergence of Algorithm \ref{ALGRPFM}.
Without loss of generality, we assume $M_{\upsilon} = M$ in the following
analysis. Namely, we do not discriminate between $M_{\upsilon}$ and $M$.
Similarly to the analysis techniques of the references
\cite{Wright1994,Xu1991,Zhang1994}, we firstly construct an auxiliary sequence
$(u^{k}, \, v^{k})$ as follows:
\begin{align}
      u^{k+1} = u^{k} + \alpha_{k} (x^{k} + \Delta x^{k} - u^{k}),
      \label{AUXSEQUK} \\
      v^{k+1} = v^{k} + \alpha_{k} (y^{k} + \Delta y^{k} - v^{k}),
      \label{AUXSEQVK}
\end{align}
where $(\Delta x^{k}, \, \Delta y^{k})$ are solved by equations
\eqref{CONNML}-\eqref{CONNMC} and $(u^{0}, \, v^{0}) \le (x^{0}, \, y^{0})$
satisfies the feasibility $r_{q}(u^{0}, \, v^{0}) = 0$. Then, it is not difficult
to verify
\begin{align}
    & r_{q}(u^{k}, \, v^{k}) = v^{k} - (Mu^{k} + q) = 0, \label{LCPFUV} \\
    & x^{k+1} - u^{k+1} = (1 - \alpha_{k}) (x^{k} - u^{k}) \ge 0, \label{XKGEUK} \\
    & y^{k+1} - v^{k+1} = (1 - \alpha_{k})(y^{k} - v^{k}) \ge 0, \label{YKGEVK}
\end{align}
where $\alpha_{k} = \Delta t_{k}/(1+\Delta t_{k})$.

\vskip 2mm

Meanwhile, in order to obtain the global convergence, we need to enforce the
condition specified below. Firstly, we select a constant $\gamma \in (0, \, 1)$
that satisfies
\begin{align}
    0 < \gamma \le \frac{\mu_{0}}
    {\min_{1 \le i \le n} \{x^{0}_{i}y^{0}_{i}\}}
    = \frac{(x^{0})^{T}y^{0} + \|r_{q}^{0}\|}
    {2n \, \min_{1 \le i \le n} \{x^{0}_{i}y^{0}_{i}\}}. \label{DEFGAMMA}
\end{align}
Then, we select $\alpha_{k}$ such that
\begin{align}
     x^{k}_{i}(\alpha)y^{k}_{i}(\alpha) \ge \gamma \mu_{k}(\alpha),
     \; i = 1, \, 2, \, \ldots, n     \label{XYGEUK}
\end{align}
holds for all $\alpha \in (0, \; \alpha_{k}] \subset (0, \, 1]$, where
$x^{k}(\alpha)$, $y^{k}(\alpha)$, $r_{q}^{k}(\alpha)$ and $\mu_{k}(\alpha)$ are
defined by
\begin{align}
    & x^{k}(\alpha) = x^{k} + \alpha \Delta x^{k}, \;
    y^{k}(\alpha) = y^{k} + \alpha \Delta y^{k}, \label{XKYKALPHA} \\
    & r_{q}^{k}(\alpha) = r_{q}(x^{k}(\alpha), y^{k}(\alpha)), \;
    \mu_{k}(\alpha) = \frac{(x^{k}(\alpha))^{T}y^{k}(\alpha) + \|r_{q}^{k}(\alpha)\|}{2n}.
    \label{RMUKALPHA}
\end{align}

\vskip 2mm

Condition \eqref{XYGEUK} is to prevent iterations from prematurely getting too
close to the boundary of the positive quadrant and its restriction on $\gamma$
is very mild. In the practice, it can be selected to be very small.

\vskip 2mm

In order to establish the main global convergence of Algorithm \ref{ALGRPFM},
similarly to the results \cite{Wright1994,Zhang1994}, we prove the following
several technique lemmas.

\vskip 2mm

\begin{lemma} \label{LEMXYBOUD}
When $(x^{k}(\alpha), \, y^{k}(\alpha))$ satisfies the condition
\eqref{XYGEUK}, where $(\Delta x^{k}, \, \Delta y^{k})$ is the solution of
equation \eqref{CONNML}-\eqref{CONNMC}, we have
\begin{align}
    (x^{k}(\alpha), \, y^{k}(\alpha)) \ge 0. \label{XYALPHAGE0}
\end{align}
Furthermore, when $x^{k}_{i}y^{k}_{i} \ge \gamma \mu_{k} \;
(i=1, \, 2, \ldots, \, n)$ and $\alpha$ satisfies
\begin{align}
    0 \le \alpha \le \min\left\{1, \; \frac{1 - \gamma/2}{1+\gamma/2}
    \frac{\sigma_{k}\mu_{k}}{\left\|\Delta X^{k}\Delta y^{k}\right\|_{\infty}}
    \right\}, \label{ALPAUBND}
\end{align}
the inequality \eqref{XYGEUK} holds.
\end{lemma}
\proof By summing two sides of the condition \eqref{XYGEUK}, we obtain
\begin{align}
    & (x^{k}(\alpha))^{T}y^{k}(\alpha) \ge n \gamma \mu_{k}(\alpha)
    = \frac{1}{2}\gamma
    \left((x^{k}(\alpha))^{T}y^{k}(\alpha)+ (1-\alpha)\|r_{q}^{k}\|\right)
    \nonumber \\
    & \hskip 2mm \ge \frac{1}{2} \gamma (x^{k}(\alpha))^{T}y^{k}(\alpha).
    \label{XYKGEGXYK}
\end{align}
Consequently, from equation \eqref{XYKGEGXYK} and $0 < \gamma < 1$, we obtain
$(x^{k}(\alpha))^{T}y^{k}(\alpha) \ge 0$. By substituting it into the condition
\eqref{XYGEUK}, we have $x^{k}_{i}(\alpha)y^{k}_{i}(\alpha) \ge 0 \,
(i = 1, \, 2, \, \ldots, \, n)$. From equation \eqref{CONNMC},
we obtain
\begin{align}
   & \alpha \sigma_{k} \mu_{k} = \alpha x^{k}_{i}y^{k}_{i}
   + \alpha x^{k}_{i} \Delta y^{k}_{i} + \alpha y^{k}_{i}\Delta x^{k}_{i}
   \nonumber \\
   & \hskip 2mm = (\alpha - 2) x^{k}_{i}y^{k}_{i} + y^{k}_{i}x^{k}_{i}(\alpha)
   + x^{k}_{i}y^{k}_{i}(\alpha), \; i = 1, \, 2, \ldots, \, n. \label{POSEQCONC}
\end{align}
If $(x^{k}_{i} (\alpha), \, y^{k}_{i}(\alpha)) < 0$ or
$x^{k}_{i}(\alpha) = 0, \, y^{k}_{i}(\alpha) < 0$ or
$x^{k}_{i}(\alpha) < 0, \, y^{k}_{i}(\alpha) = 0$, by combining it with
$(x^{k}_{i}, \, y^{k}_{i}) \ge 0$, it contradicts equation
\eqref{POSEQCONC}. Therefore, we obtain $(x^{k}_{i}(\alpha), \, y^{k}_{i}(\alpha))
\ge 0 \, (i = 1, \, 2, \, \ldots, \, n)$, which proves equation \eqref{XYALPHAGE0}.

\vskip 2mm

From equations \eqref{CONNML}-\eqref{CONNMC}, we have
\begin{align}
    &x^{k}_{i}(\alpha)y^{k}_{i}(\alpha) - \gamma \mu_{k}(\alpha)
    = x_{i}^{k}y_{i}^{k} +
    \alpha(x^{k}_{i} \Delta y^{k}_{i} + y^{k}_{i} \Delta x^{k}_{i})
    + \alpha^{2}\Delta x^{k}_{i} \Delta y^{k}_{i} \nonumber \\
    & \hskip 2mm - \gamma \mu_{k} - \frac{1}{2n}\gamma \alpha
    \left((x^{k})^{T}\Delta y^{k} + (y^{k})^{T} \Delta x^{k} - \|r_{q}^{k}\|
    + \alpha (\Delta x^{k})^{T} \Delta y^{k}\right) \nonumber \\
    & = x_{i}^{k}y_{i}^{k} +
    \alpha(\sigma_{k} \mu_{k} - x^{k}_{i} y^{k}_{i})
    + \alpha^{2}\Delta x^{k}_{i} \Delta y^{k}_{i} \nonumber \\
    & \hskip 2mm - \gamma \mu_{k} - \frac{1}{2n}\gamma \alpha
    \left(n \sigma_{k} \mu_{k} - (x^{k})^{T}y^{k} - \|r_{q}^{k}\|
    + \alpha (\Delta x^{k})^{T} \Delta y^{k}\right) \nonumber \\
    & = (1-\alpha)\left(x^{k}_{i}y^{k}_{i} - \gamma \mu_{k}\right)
    + \alpha \left(\sigma_{k} \mu_{k} \left(1 - \frac{\gamma}{2}\right)
    + \alpha \left(\Delta x^{k}_{i} \Delta y^{k}_{i}
    - \frac{\gamma}{2n}(\Delta x^{k})^{T}\Delta y^{k}\right)\right) \nonumber \\
    & \ge (1-\alpha)\left(x^{k}_{i}y^{k}_{i} - \gamma \mu_{k}\right)
    + \alpha \left(\sigma_{k} \mu_{k} \left(1 - \frac{\gamma}{2}\right)
    - \alpha \left(1 + \frac{\gamma}{2}\right)
    \left\|\Delta X^{k} \Delta y^{k} \right\|_{\infty}\right), \label{WNNEQALP}
\end{align}
where $\Delta X^{k} = \text{diag}(\Delta x^{k})$. By substituting
$x^{k}_{i}y^{k}_{i} \ge \gamma \mu_{k}$ and $0 < \alpha \le 1$ into
equation \eqref{WNNEQALP}, we have
\begin{align}
    x^{k}_{i}(\alpha)y^{k}_{i}(\alpha) - \gamma \mu_{k}(\alpha)
    \ge \alpha \left(\sigma_{k} \mu_{k} \left(1 - \frac{\gamma}{2}\right)
    - \alpha \left(1 + \frac{\gamma}{2}\right)
    \left\|\Delta X^{k} \Delta y^{k} \right\|_{\infty}\right).
    \label{XYGEUALP}
\end{align}
Therefore, when $\alpha$ satisfies equation \eqref{ALPAUBND},
from equation \eqref{XYGEUALP}, we obtain
\begin{align}
    x^{k}_{i}(\alpha)y^{k}_{i}(\alpha) - \gamma \mu_{k}(\alpha) \ge 0, \;
    i = 1, \, 2, \, \ldots, \, n, \nonumber
\end{align}
which gives the inequality \eqref{XYGEUK}. \eproof

\vskip 2mm

\begin{lemma} \label{LEMXYKUB}
Assume that $\{(x^{k}, \, y^{k})\}$ is generated by Algorithm \ref{ALGRPFM}.
Then, for any feasible solution $(\bar{x}, \, \bar{y})$ of the linear
complementarity problem \eqref{NLELCP}, there exists a positive constant
$C_{1}$ such that
\begin{align}
   & (y^{k})^{T}(x^{k} - u^{k}) + (x^{k})^{T}(y^{k} - v^{k})
   \le C_{1}  \label{YXMULEUB}
\end{align}
holds for all $k = 0, \, 1, \, \ldots$. Furthermore, if $(\bar{x}, \, \bar{y})$
is strictly feasible, $\{(x^{k}, \, y^{k})\}$ is bounded.
\end{lemma}
\proof From equation \eqref{LCPFUV}, we know that $(u^{k}, \, v^{k})$ satisfies
the feasibility. By combining it with the semi-definite positivity of $M$, we
have
\begin{align}
     (\bar{x} - u^{k})^{T}(\bar{y} - v^{k})
     = (\bar{x} - u^{k})^{T}M(\bar{x} - u^{k}) \ge 0. \label{DXUYVGEZ}
\end{align}
Consequently, from equation \eqref{DXUYVGEZ}, we have
\begin{align}
    0 & \le (\bar{x} - u^{k})^{T}(\bar{y} - v^{k}) \nonumber \\
    & = (\bar{x} - x^{k}+ (x^{k} - u^{k}))^{T}(\bar{y} - y^{k} + (y^{k} - v^{k}))
    \nonumber \\
    &
    = \bar{x}^{T}(y^{k} - v^{k})  - (x^{k})^{T}(y^{k} - v^{k})
    + \bar{x}^{T}\bar{y} + (x^{k})^{T}y^{k} - \bar{x}^{T}y^{k} - \bar{y}^{T}x^{k}
    \nonumber \\
    & \hskip 2mm
    + (x^{k} - u^{k})^{T}(y^{k} - v^{k}) + \bar{y}^{T}(x^{k} - u^{k})
    - (y^{k})^{T}(x^{k} - u^{k}) . \label{SPLGEZ}
\end{align}
From equations \eqref{XKGEUK}-\eqref{YKGEVK}, we have
\begin{align}
    0 \le x^{k}- u^{k} \le x^{0} - u^{0}, \; 0 \le y^{k}- v^{k} \le y^{0} - v^{0}.
    \label{DXULEUPB}
\end{align}
By substituting equation \eqref{DXULEUPB} into equation \eqref{SPLGEZ}, from
$(x^{k}, \, y^{k}) \ge 0, \; (\bar{x}, \, \bar{y}) \ge 0$ and
$\phi(x^{k}, \, y^{k}) \le \phi(x^{0}, \, y^{0})$, we obtain
\begin{align}
    & (y^{k})^{T}(x^{k} - u^{k})  + (x^{k})^{T}(y^{k} - v^{k})
    \le \bar{x}^{T}(y^{k} - v^{k}) + \bar{y}^{T}(x^{k} - u^{k}) \nonumber \\
    & \hskip 2mm
    + (x^{k})^{T}y^{k} + (x^{k} - u^{k})^{T}(y^{k} - v^{k}) + \bar{x}^{T}\bar{y}
    \nonumber \\
    & \hskip 2mm
    \le \bar{x}^{T}(y^{0} - v^{0}) + \bar{y}^{T}(x^{0} - u^{0})
    + \phi(x^{k}, \, y^{k}) + (x^{0} - u^{0})^{T}(y^{0} - v^{0})
    + \bar{x}^{T}\bar{y}  \nonumber \\
    & \hskip 2mm
    \le \bar{x}^{T}(y^{0} - v^{0}) + \bar{y}^{T}(x^{0} - u^{0})
    + \phi(x^{0}, \, y^{0}) + (x^{0} - u^{0})^{T}(y^{0} - v^{0})
    + \bar{x}^{T}\bar{y} \nonumber \\
    & \hskip 2mm
    = \bar{x}^{T}(y^{0} - v^{0}) + \bar{y}^{T}(x^{0} - u^{0})
    + (x^{0} - u^{0})^{T}(y^{0} - v^{0}) + 2n\mu_{0}
    + \bar{x}^{T}\bar{y}. \label{YDXULEUPB}
\end{align}
We set
\begin{align}
   C_{1} = (y^{0} - v^{0})^{T}(x^{0} - u^{0}) + \bar{y}^{T}(x^{0}
    - u^{0}) + \bar{x}^{T}(y^{0} - v^{0}) + \bar{x}^{T}\bar{y}+ 2n\mu_{0}.
    \nonumber
\end{align}
Then, from equation \eqref{YDXULEUPB}, we obtain the inequality \eqref{YXMULEUB}.

\vskip 2mm

When $(\bar{x}, \, \bar{y})$ is strictly feasible, from inequalities
\eqref{SPLGEZ}-\eqref{DXULEUPB}, we have
\begin{align}
    & \bar{y}^{T}x^{k} + \bar{x}^{T}y^{k} \le \bar{x}^{T}(y^{k} - v^{k})
    + \bar{y}^{T}(x^{k} - u^{k}) + \bar{x}^{T}\bar{y} + (x^{k})^{T}y^{k}
    + (x^{k} - u^{k})^{T}(y^{k} - v^{k}) \nonumber \\
    & \hskip 2mm
    \le \bar{x}^{T}(y^{0} - v^{0})
    + \bar{y}^{T}(x^{0} - u^{0}) + \bar{x}^{T}\bar{y} + \phi(x^{k}, \, y^{k})
    + (x^{0} - u^{0})^{T}(y^{0} - v^{0}) \nonumber \\
    & \hskip 2mm
    \le \bar{x}^{T}(y^{0} - v^{0})
    + \bar{y}^{T}(x^{0} - u^{0}) + \bar{x}^{T}\bar{y} + 2n \mu_{0}
    + (x^{0} - u^{0})^{T}(y^{0} - v^{0}) \triangleq C_{2}, \label{XYKLEUPB}
\end{align}
where we use the property $\phi(x^{k}, \, y^{k}) \le \phi(x^{0}, \, y^{0})
= 2n\mu_{0}$ in the third inequality. Consequently, from equation
\eqref{XYKLEUPB} and $(\bar{x}, \, \bar{y}) > 0$, we have
\begin{align}
    0 \le x^{k}_{i} \le \frac{C_{2}}{\min_{1 \le i \le n} \{\bar{y}_{i}\}}, \;
    0 \le y^{k}_{i} \le \frac{C_2}{\min_{1 \le i \le n} \{\bar{x}_{i}\}},
    \; i = 1, \, 2, \, \ldots, \, n, \nonumber
\end{align}
which prove the boundedness of $(x^{k}, \, y^{k})$. \eproof

\vskip 2mm

\begin{lemma} \label{LEMDXYLEUB}
Assume that $\{(x^{k}, \, y^{k})\}$ is generated by Algorithm \ref{ALGRPFM} and
satisfies the condition \eqref{XYGEUK}. Then, when $\mu_{k} \ge \epsilon > 0$,
there exists a positive constant $\omega^{\ast}$ such that
\begin{align}
   \|D^{k}\Delta x^{k}\|^{2} + \|(D^{k})^{-1}\Delta y^{k}\|^{2}
   \le (\omega^{\ast})^{2},    \label{DXYLEUB}
\end{align}
holds for all $k = 0, \, 1, \, \ldots$, where $D^{k} =
\text{diag}\left((y^{k}./x^{k})^{1/2}\right)$, $(D^{k})^{-1} =
\text{diag}\left((x^{k}./y^{k})^{1/2}\right)$ and
$(\Delta x^{k}, \, \Delta y^{k})$ is the solution of equations
\eqref{CONNML}-\eqref{CONNMC}.
\end{lemma}

\proof We denote
\begin{align}
   \omega_{k} = \left( \|D^{k}\Delta x^{k}\|^{2}
   + \|(D^{k})^{-1}\Delta y^{k}\|^{2}\right)^{1/2}. \label{DEFTK}
\end{align}
Then, from equation \eqref{DEFTK}, we have
\begin{align}
    \|D^{k}\Delta x^{k}\|_{\infty} \le \|D^{k}\Delta x^{k}\| \le \omega_{k}, \;
    \|(D^{k})^{-1}\Delta y^{k}\|_{\infty} \le \|(D^{k})^{-1}\Delta y^{k}\|
    \le \omega_{k}. \label{DKDELXKOM}
\end{align}

\vskip 2mm

From equation \eqref{LCPFUV}, we know that $(u^{k}, \, v^{k})$ is feasible. By
combining it with equation \eqref{CONNML}, we obtain
\begin{align}
     y^{k} + \Delta y^{k} - v^{k} = M(x^{k} + \Delta x^{k} - u^{k}).
     \label{YK1MVK}
\end{align}
Therefore, from equation \eqref{YK1MVK} and the semi-definite positivity of $M$,
we have
\begin{align}
    (y^{k} + \Delta y^{k} - v^{k})^{T}(x^{k} + \Delta x^{k} - u^{k})
    = (x^{k} + \Delta x^{k} - u^{k})^{T}M(x^{k} + \Delta x^{k} - u^{k})
    \ge 0. \label{YDLEYMYGEZ}
\end{align}
From equations \eqref{XKGEUK}-\eqref{YKGEVK}, we know
\begin{align}
     x^{0} - u^{0} \ge x^{k} - u^{k}, \; y^{0} - v^{0} \ge y^{k} - v^{k}. \nonumber
\end{align}
By substituting them into equation \eqref{YDLEYMYGEZ}, we obtain
\begin{align}
   & (x^{0} - u^{0})^{T}(y^{0} - v^{0}) + (y^{k} - v^{k})^{T}\Delta x^{k}
   + (x^{k} - u^{k})^{T} \Delta y^{k} + (\Delta x^{k})^{T}\Delta y^{k}
   \nonumber \\
   & \ge
   (x^{k} - u^{k})^{T}(y^{k} - v^{k}) + (y^{k} - v^{k})^{T}\Delta x^{k}
   + (x^{k} - u^{k})^{T} \Delta y^{k} + (\Delta x^{k})^{T}\Delta y^{k}
   \ge 0. \label{XMUDXYGEZ}
\end{align}

\vskip 2mm

By using the inequality $|x^{T}y| \le \|x\| \|y\| \le \|x\|_{1}\|y\|$,
from equation \eqref{DKDELXKOM}, we have
\begin{align}
    & \left|(y^{k} - v^{k})^{T}\Delta x^{k}\right| =
    \left|\left((D^{k})^{-1}(y^{k} - v^{k})\right)^{T}(D^{k}\Delta x^{k})\right| \nonumber \\
    & \hskip 2mm
    \le \left\|(D^{k})^{-1}(y^{k} - v^{k})\right\|_{1}\|D^{k}\Delta x^{k}\|
    \le \left\|(D^{k})^{-1}(y^{k} - v^{k})\right\|_{1} \omega_{k}.
    \label{YMVDELXLEB}
\end{align}
From the condition \eqref{XYGEUK} and $y^{k} - v^{k} \ge 0$, we have
\begin{align}
   ({x^{k}_{i}}/{y^{k}_{i}})^{1/2}(y^{k}_{i} - v^{k}_{i}) =
   x^{k}_{i}(y^{k}_{i} - v^{k}_{i})/(x^{k}_{i}y^{k}_{i})^{1/2}
   \le x^{k}_{i}(y^{k}_{i} - v^{k}_{i})/(\gamma \mu_{k})^{1/2}. \label{XKMXMVDY}
\end{align}
By substituting equation \eqref{XKMXMVDY} into equation \eqref{YMVDELXLEB}, we
obtain
\begin{align}
    & \left|(y^{k} - v^{k})^{T}\Delta x^{k}\right| \le
    \left(\sum_{i=1}^{n}({x^{k}_{i}}/{y^{k}_{i}})^{1/2}(y^{k}_{i} - v^{k}_{i})\right)
    \omega_{k} \nonumber \\
    & \hskip 2mm
    \le (x^{k})^{T}(y^{k} - v^{k})\omega_{k}/(\gamma \mu_{k})^{1/2}. \label{YMVDELUB}
\end{align}
Similarly to the proof of equation \eqref{YMVDELUB}, we obtain
\begin{align}
    \left|(x^{k} - u^{k})^{T}\Delta y^{k}\right| \le
     (y^{k})^{T}(x^{k} - u^{k})\omega_{k}/(\gamma \mu_{k})^{1/2}. \label{XMUDELUB}
\end{align}
By substituting inequalities \eqref{YXMULEUB} and \eqref{YMVDELUB}-\eqref{XMUDELUB}
into inequality \eqref{XMUDXYGEZ}, we obtain
\begin{align}
   & (\Delta x^{k})^{T}\Delta y^{k} \ge - (x^{0} - u^{0})^{T}(y^{0} - v^{0})
   - \left|(y^{k} - v^{k})^{T}\Delta x^{k}\right|
   - \left|(x^{k} - u^{k})^{T} \Delta y^{k}\right| \nonumber \\
   & \hskip 2mm
   \ge - (x^{0} - u^{0})^{T}(y^{0} - v^{0}) - ((y^{k})^{T}(x^{k} - u^{k})
   + (x^{k})(y^{k} - v^{k}))\omega_{k}/(\gamma \mu_{k})^{1/2} \nonumber \\
  & \hskip 2mm
  \ge - (x^{0} - u^{0})^{T}(y^{0} - v^{0})
  - \left(C_{1}/(\gamma \epsilon)^{1/2}\right)\omega_{k}. \label{DELXYGEN}
\end{align}

\vskip 2mm

From equation \eqref{CONNMC} and the condition \eqref{XYGEUK}, we have
\begin{align}
  & \|D^{k}\Delta x^{k}\|^{2} + \|(D^{k})^{-1} \Delta y^{k}\|^{2}
  + 2 (\Delta x^{k})^{T} \Delta y^{k}
  = \|D^{k}\Delta x^{k} + (D^{k})^{-1} \Delta y^{k}\|^{2} \nonumber \\
  & \hskip 2mm
  = (\sigma_{k}\mu_{k})^{2} \sum_{i=1}^{n}\frac{1}{x^{k}_{i}y^{k}_{i}}
  + (x^{k})^{T}y^{k} - 2\sigma_{k}\mu_{k} n \nonumber \\
  & \hskip 2mm
  \le n \sigma_{k}^{2} \mu_{k} /\gamma
  + \phi(x^{k}, \, y^{k}) - 2\sigma_{k}\mu_{k} n
  \le n \mu_{0} \left(\sigma_{max}^{2}/\gamma + 2 - 2\sigma_{min}\right).
  \label{DELXYLEMUK}
\end{align}
By substituting inequality \eqref{DELXYGEN}  into inequality \eqref{DELXYLEMUK},
we obtain
\begin{align}
   q(\omega_{k}) \triangleq \omega_{k}^{2}
   - 2\left(C_{1}/(\gamma\epsilon)^{1/2}\right)\omega_{k} - \zeta \le 0,
   \label{QUAFUN}
\end{align}
where the positive constant $\zeta$ is defined by
\begin{align}
    \zeta = 2(x^{0} - u^{0})^{T}(y^{0} - v^{0})
    + n \mu_{0}(\sigma_{max}^{2}/\gamma + 2 - 2\sigma_{min}). \label{DEFZETA}
\end{align}
The quadratic function $q(\omega)$ is convex and has a unique positive root at
\begin{align}
   \omega^{\ast} = C_{1}/(\gamma \epsilon)^{1/2}
   + \sqrt{C_{1}^{2}/(\gamma \epsilon) + \zeta}. \nonumber
\end{align}
This implies
\begin{align}
    \omega_{k} \le \omega^{\ast}. \label{OMEGAKLEUB}
\end{align}
Consequently, we obtain inequality \eqref{DXYLEUB}. \eproof

\vskip 2mm

In order to prove the global convergence of Algorithm \ref{ALGRPFM}, we need to
estimate the positive lower bound of $\Delta t_{k} \, (k = 0, \, 1, \, \ldots)$.

\begin{lemma} \label{LEMDTKGEZ}
Assume that $\{(x^{k}, \, y^{k})\}$ is generated by Algorithm \ref{ALGRPFM} and
satisfies the condition \eqref{XYGEUK}. Then, when $\mu_{k} \ge \epsilon > 0$,
there exists a positive constant $\delta_{\Delta t}$ such that
\begin{align}
     \Delta t_{k} \ge \frac{1}{2}\delta_{\Delta t} > 0   \label{DELTATKGEZ}
\end{align}
holds for all $k = 0, \, 1, \, \ldots$.
\end{lemma}
\proof From inequality \eqref{DXYLEUB}, we have
\begin{align}
    & |\Delta x^{k}_{i} \Delta y^{k}_{i}|
    = \left|(y^{k}_{i}/x^{k}_{i})^{1/2}\Delta x^{k}_{i}\right|
    \left|(x^{k}_{i}/y^{k}_{i})^{1/2} \Delta y^{k}_{i}\right|
    \le \|D^{k}\Delta x^{k}\| \|(D^{k})^{-1} \Delta y^{k}\| \nonumber \\
    & \hskip 2mm
    \le \frac{1}{2}\omega_{k}^{2} \le \frac{1}{2} (\omega^{\ast})^{2}, \; i = 1, \, 2, \, \ldots, n,
    \nonumber
\end{align}
which gives
\begin{align}
    \|\Delta X^{k} \Delta y^{k}\|_{\infty} \le \frac{1}{2} (\omega^{\ast})^{2}, \;
    k = 0, \, 1, \, 2, \, \ldots.     \label{DELXYLEOME}
\end{align}
Consequently, we have
\begin{align}
     \frac{\sigma_{k}\mu_{k}}{\|\Delta X^{k} \Delta y^{k}\|_{\infty}}
     \ge \frac{2\sigma_{min} \, \epsilon}{(\omega^{\ast})^{2}}, \;
     k = 0, \, 1, \, 2, \, \ldots. \label{SIGMUKGEUB}
\end{align}
We denote
\begin{align}
   \alpha^{\ast}_{\mu} = \min \left\{1, \; \frac{1 - \gamma/2}{1 + \gamma/2}
   \, \frac{(2\sigma_{min} \, \epsilon)}{(\omega^{\ast})^{2}} \right\}. \label{ALPHAAST}
\end{align}
Then, when $\Delta t_{k}$ satisfies
\begin{align}
     \Delta t_{k} \le \frac{\alpha^{\ast}_{\mu}}{1 - \alpha^{\ast}_{\mu}},
     \label{DELTKLEALP}
\end{align}
from equations \eqref{SIGMUKGEUB}-\eqref{ALPHAAST} and \eqref{ALPAUBND}, we know
that the condition \eqref{XYGEUK} holds, where
$\alpha_{k} = \Delta t_{k}/(1 + \Delta t_{k})$.

\vskip 2mm

From equation \eqref{DXYLEUB} and the Cauchy-Schwartz inequality
$|x^{T}y| \le \|x\| \|y\|$, we have
\begin{align}
   & |(\Delta x^{k})^{T}\Delta y^{k}|
   = \left|(D^{k}\Delta x^{k})^{T}((D^{k})^{-1}\Delta y^{k})\right|
   \le \|D^{k}\Delta x^{k}\| \|(D^{k})^{-1}\Delta y^{k}\| \nonumber \\
   & \hskip 2mm
   \le \frac{1}{2}\left(\|D^{k}\Delta x^{k}\|^{2}
   + \|(D^{k})^{-1}\Delta y^{k}\|^{2}\right)
   \le \frac{1}{2} (\omega^{\ast})^{2}. \label{DELXYDLEUB}
\end{align}
From equation \eqref{CONNMC}, we have
\begin{align}
   (x^{k})^{T}\Delta y^{k} + (y^{k})^{T}\Delta x^{k}
   = n\sigma_{k}\mu_{k} - (x^{k})^{T}y^{k}. \label{XDYAYDX}
\end{align}
By substituting equations \eqref{DELXYDLEUB}-\eqref{XDYAYDX} into equation
\eqref{RHOK}, we obtain
\begin{align}
   & \rho_{k} = \frac{Ared_{k}}{Pred_{k}} =
   1 - \frac{\alpha_{k}(\Delta x^{k})^{T}\Delta y^{k}}
   {\|r_{q}^{k}\| + (x^{k})^{T}y^{k}  - n\sigma_{k}\mu_{k}} \nonumber \\
   & \hskip 2mm
   = 1 - \frac{\alpha_{k}(\Delta x^{k})^{T}\Delta y^{k}}
   {(2 - \sigma_{k})n \mu_{k}}
   \ge 1 - \frac{(\omega^{\ast})^{2}}{2(2-\sigma_{max}) n \epsilon} \alpha_{k}.
   \label{RHOKLEUB}
\end{align}
We denote
\begin{align}
   \alpha_{\rho}^{\ast} = \min \left\{1, \;
   \frac{2(2 - \sigma_{max}) (1-\eta_{2})n \epsilon}
   {(\omega^{\ast})^{2}}\right\}. \label{ALPHARA}
\end{align}
Then, when $\Delta t_{k}$ satisfies
\begin{align}
    \Delta t_{k} \le \frac{\alpha_{\rho}^{\ast}}
    {1 - \alpha_{\rho}^{\ast}}, \label{DTLEALPR}
\end{align}
from equations \eqref{RHOKLEUB}-\eqref{ALPHARA}, we know $\rho_{k} \ge \eta_2$,
where $\alpha_{k} = \Delta t_{k}/(1 + \Delta t_{k})$.

\vskip 2mm

We denote
\begin{align}
    \delta_{\Delta t} \triangleq \min\left\{\frac{\sigma_{\mu}^{\ast}}
    {1 - \sigma_{\mu}^{\ast}}, \; \frac{\alpha_{\rho}^{\ast}}
    {1 - \alpha_{\rho}^{\ast}}, \; \Delta t_{0}\right\}.   \label{PARDELAT}
\end{align}
Assume that $K$ is the first index such that $\Delta t_{K}
\le \delta_{\Delta t}$ holds.
Then, from equations \eqref{DELTKLEALP}-\eqref{PARDELAT}, we know that
the condition \eqref{XYGEUK} holds and $\rho_{K} \ge \eta_{2}$.
Consequently, according to the time-stepping adjustment scheme \eqref{TSK1},
$\Delta t_{K+1}$ will be enlarged. Therefore,
$\Delta t_k \ge (1/2) \delta_{\Delta t}$ holds for all
$k = 0, \, 1, \, 2, \, \dots$. \qed

\vskip 2mm

By using the estimation results of Lemma \ref{LEMDTKGEZ}, we prove
that $\{\mu_{k}\}$ converges to zero when $k$ tends to infinity as follows.

\begin{theorem} \label{THECOVFUK}
Assume that $\{(x^{k}, \, y^{k})\}$ is generated by Algorithm \ref{ALGRPFM} and
satisfies the condition \eqref{XYGEUK}. Then, we have
\begin{align}
    \lim_{k \to \infty} \mu_{k} = 0. \label{MUKTOZ}
\end{align}
\end{theorem}
\proof We prove the result \eqref{MUKTOZ} by contradiction.
Assume that there exists a positive constant $\epsilon$ such that
\begin{align}
     \mu_{k} \ge \epsilon > 0 \label{MUKGEZ}
\end{align}
holds for all $k = 0, \, 1, \, \ldots$. Then, according to Algorithm \ref{ALGRPFM}
and Lemma \ref{LEMDTKGEZ}, we know that
there exists an infinite subsequence $\{(x^{k_l}, \, y^{k_l})\}$ such that
\begin{align}
   \frac{\phi(x^{k_l}, \, y^{k_l}) - \phi(x^{k_l} + \alpha_{k_l} \Delta x^{k_l},
   \, y^{k_l} + \alpha_{k_l} \Delta y^{k_l})}
   {\phi(x^{k_l}, \, y^{k_l}) - m_{k_l}(\alpha_{k_l})} \ge \eta_{1}
   \label{RHOKGEPC}
\end{align}
holds for all $l = 0, \, 1, \, 2, \, \ldots$, where $\alpha_{k_l} =
\Delta t_{k_l}/(1 + \Delta t_{k_l})$. Otherwise, all steps are rejected
after a given iteration index, then the time step  will keep decreasing,
which contradicts equation \eqref{DELTATKGEZ}.

\vskip 2mm

From equations \eqref{DELTATKGEZ}, \eqref{RHOKLEUB} and \eqref{RHOKGEPC}, we have
\begin{align}
    & \phi(x^{k_l}, \, y^{k_l}) - \phi(x^{k_l} + \alpha_{k_l} \Delta x^{k_l},
   \, y^{k_l} + \alpha_{k_l} \Delta y^{k_l}) \nonumber \\
    & \hskip 2mm
    \ge \frac{\eta_{1}\Delta t_{k_l}}{1+\Delta t_{k_l}}(2 - \sigma_{k_l})n \mu_{k_l}
    \ge  \frac{\eta_{1}\delta_{\Delta t}/2}{1+ \delta_{\Delta t}/2}
    (2 - \sigma_{max})n \mu_{k_l}.    \label{PHIKDGE}
\end{align}
Therefore, from equation \eqref{PHIKDGE} and $\phi(x^{k+1}, \, y^{k+1}) \le
\phi(x^{k}, \, y^{k})$, we have
\begin{align}
    & \phi(x^{0}, \, y^{0}) \ge \phi(x^{0}, \, y^{0}) - \lim_{k \to \infty}
    \phi(x^{k}, \, y^{k}) =
    \sum_{k = 0}^{\infty} \left(\phi(x^{k}, \, y^{k}) - \phi(x^{k+1}, \, y^{k+1})\right)
    \nonumber \\
    & \hskip 2mm \ge \sum_{l = 0}^{\infty} \left(\phi(x^{k_l}, \, y^{k_l})
    - \phi(x^{k_l + 1}, \, y^{k_l + 1})\right)
    \ge  \frac{\eta_{1} \delta_{\Delta t}/2}{1+\delta_{\Delta t}/2}
    \sum_{l=0}^{\infty}(2 - \sigma_{max})n \mu_{k_l}.
  \label{PHI0GESUM}
\end{align}
Consequently, from equation \eqref{PHI0GESUM}, we obtain
\begin{align}
    \lim_{l \to \infty} \mu_{k_l} = 0,  \nonumber
\end{align}
which contradicts the assumption $\mu_{k} \ge \epsilon > 0$ for all
$k = 0, \, 1, \ldots$. Therefore, we have
\begin{align}
    \lim_{k \to \infty} \inf \mu_{k} = 0. \label{INFMUKTOZ}
\end{align}
Since $\mu_{k} = \phi(x^{k}, \, y^{k})/(2n)$ and $\{\phi(x^{k}, \, y^{k})\}$
is monotonically decreasing, we know that $\{\mu_{k}\}$ is monotonically decreasing.
By combining it with equation \eqref{INFMUKTOZ}, we know that the result
\eqref{MUKTOZ} is true. \qed

\begin{remark}
In the analysis framework of the global convergence, in comparison to
that of the known path-following algorithm, we do not need to select
$\alpha_{k}$ such that the condition of the priority to feasibility over
complementarity
\begin{align}
    \|r_{q}^{k}(\alpha)\| \le (x^{k}(\alpha))^{T}y^{k}(\alpha)
    \frac{\|r_{q}^{0}\|}{(x^{0})^{T}y^{0}} \label{FOVERC}
\end{align}
holds for all $\alpha \in (0, \; \alpha_{k}] \subset (0, \; 1]$.
\end{remark}

\vskip 2mm

\begin{theorem} \label{THEXYKGCONV}
Assume that $\{(x^{k}, \, y^{k})\}$ is the sequence of iterations generated by
Algorithm \ref{ALGRPFM} and satisfies the condition \eqref{XYGEUK}. If the linear
complementarity problem has a strictly feasible solution $(\bar{x}, \, \bar{y})$,
$\{(x^{k}, \, y^{k})\}$ exists a limit point $(x^{\ast}, \, y^{\ast})$ and this
limit point satisfies the linear complementarity equation \eqref{LCP}.
\end{theorem}
\proof Since the linear complementarity problem has a strictly feasible solution,
from Lemma \ref{LEMXYKUB}, we know that $\{(x^{k}, \, y^{k})\}$ is bounded.
Consequently, the sequence $\{(x^{k}, \, y^{k})\}$ exists a convergence subsequence
$\{(x^{k_l}, \, y^{k_l})\}$ and we denote its limit point as $(x^{\ast}, \, y^{\ast})$.
By combining it with Theorem \ref{THECOVFUK}, we obtain
\begin{align}
    & \lim_{l \to \infty} \mu_{k_l} =
    \frac{1}{2n} \left(\lim_{l \to \infty} (x^{k_l})^{T}y^{k_l} + \lim_{l \to \infty}
    \|y^{k_l} - (Mx^{k_l} + q)\|\right)  \nonumber \\
    & \hskip 2mm
    = \frac{1}{2n} \left((x^{\ast})^{T}y^{\ast} + \|y^{\ast} - (Mx^{\ast} + q)\|\right)
    = 0. \label{SUBCONMUK}
\end{align}
Since $(x^{k_l}, \, y^{k_l}) \ge 0$, we have $(x^{\ast}, \, y^{\ast}) \ge 0$. By
substituting it into equation \eqref{SUBCONMUK}, we conclude
\begin{align}
    y^{\ast} - (Mx^{\ast} + q) = 0, \; (x^{\ast})^{T}y^{\ast} = 0, \;
    (x^{\ast}, \, y^{\ast}) \ge 0. \nonumber
\end{align}
It implies that $(x^{\ast}, \, y^{\ast})$ is the solution of the linear complementarity
problem \eqref{LCP}. \eproof

\section{Numerical experiments}

\vskip 2mm

In this section, some numerical experiments are conducted to test the performance
of RPFMTr (Algorithm \ref{ALGRPFM}) for the linear complementarity problems (LCPs),
in comparison to two state-of-the-art commercial solvers, i.e.
PATH \cite{DF1995,PATH,FM2000,FM2022} and MILES (a Mixed Inequality and nonLinear
Equation Solver) \cite{Mathiesen1985,Rutherford1995,Rutherford2022}. RPFMTr
(Algorithm \ref{ALGRPFM}) is coded with the MATLAB language and executed in
MATLAB (R2020a) environment \cite{MATLAB}. PATH and MILES are executed in the
GAMS v28.2 (2019) environment \cite{GAMS}. All numerical experiments are
performed by a HP notebook with the Intel quadcore CPU and 8Gb memory. MILES
solves the LCP based on the Lemke's almost complementary pivoting algorithm
\cite{CPS2009,LH1964,Rutherford1995}. PATH is a solver based on the
path-following procedure \cite{Dirkse1994,DF1995}, the Fisher's non-smooth regularization
technique \cite{Fisher1992} and the Lemke's pivoting algorithm \cite{CPS2009,LH1964}.
PATH and MILES are two robust and efficient solvers for the complementarity
problems and used in many general modelling systems such as AMPL (A Mathematical
Programming Language) \cite{AMPL} and GAMS (General Algebraic Modelling System)
\cite{GAMS}. Therefore, we select these two solvers as the basis for comparison.

\vskip 2mm

The construction approach of test problems is described as follows. Firstly, we generate a
test matrix $M$ with the following structure:
\begin{align}
    M  = \begin{bmatrix} 0 & -A^{T} \\
                         A & 0
         \end{bmatrix}, \label{TESTMAT}
\end{align}
where the test matrix $A$ comes from the linear programming subset of NETLIB
\cite{NETLIB}. It is easy to verify that the test matrix $M$ defined by equation
\eqref{TESTMAT} is semi-definite positive. Then, we generate two complementarity
vectors $x$ and $y$ as follows:
\begin{align}
    x = [1 \; 0 \; 1 \; 0 \; \cdots \; 1 \; 0]^{T}, \;
    y = [0 \; 1 \; 0 \; 1 \; \cdots \; 0 \; 1]^{T}.   \label{CVECTORS}
\end{align}
Finally, by using these two vectors $x, \, y$ defined by equation \eqref{CVECTORS}
and the test matrix $M$ defined by equation \eqref{TESTMAT}, we generate the
following test vector $q$:
\begin{align}
    q = y - Mx. \label{TESTVEC}
\end{align}
The scales of test problems vary from dimension 78 to 40216. And the termination
conditions of Algorithm \ref{ALGRPFM} (RPFMTr), PATH and MILES for finding a solution
of the nonlinear system \eqref{NLELCP} are
\begin{align}
    \|y - (Mx + q)\|_{\infty} \leq 10^{-6}, \; \|Xy\|_{\infty} \leq 10^{-6}, \;
    (x, \, y) \ge 0, \nonumber
\end{align}
where $X = \text{diag}(x)$.

\vskip 2mm

According to the manual of GAMS \cite{FM2022}, an LCP will be generated a list
of all variables appearing in the equations found in the model statement, and
the number of equations equals the number of variables. In other words, for an
LCP solved in the GAMS environment, its matrix $M$ should not include the row
with all zeros. In order to avoid this error, we add $\epsilon = 1.0 \times 10^{-6}$
to the first element of the row or column with all zeros in the test
sub-matrix $A$ of matrix $M$ defined by equation \eqref{TESTMAT}.

\vskip 2mm

Most of matrices $A$ coming from the linear programming subset of
NETLIB \cite{NETLIB} are sparse. In order to test the performance of RPFMTr,
PATH and MILES for the dense LCP further, we add a small random perturbation
to the elements of the sub-matrix $A$ in equation \eqref{TESTMAT} and generate
the dense test matrix $\overline{M}$ to replace the test matrix $M$ in equation
\eqref{TESTMAT} as follows:
\begin{align}
    \overline{M} = \begin{bmatrix}
      0 & -A_{\epsilon}^{T} \\
      A_{\epsilon} & 0
    \end{bmatrix}, \;
    A_{\epsilon} = A + \text{rand}(m_{A},\, n_{A})*\epsilon, \;
    [m_{A}, \, n_{A}] = \text{sizes}(A), \label{TESTMATD}
\end{align}
where $\epsilon = 10^{-3}$ and matrix $A$ comes from the linear programming subset
of NETLIB \cite{NETLIB}.

\vskip 2mm

\begin{figure}[!htbp]
      \centering
        \includegraphics[width=0.80\textwidth]{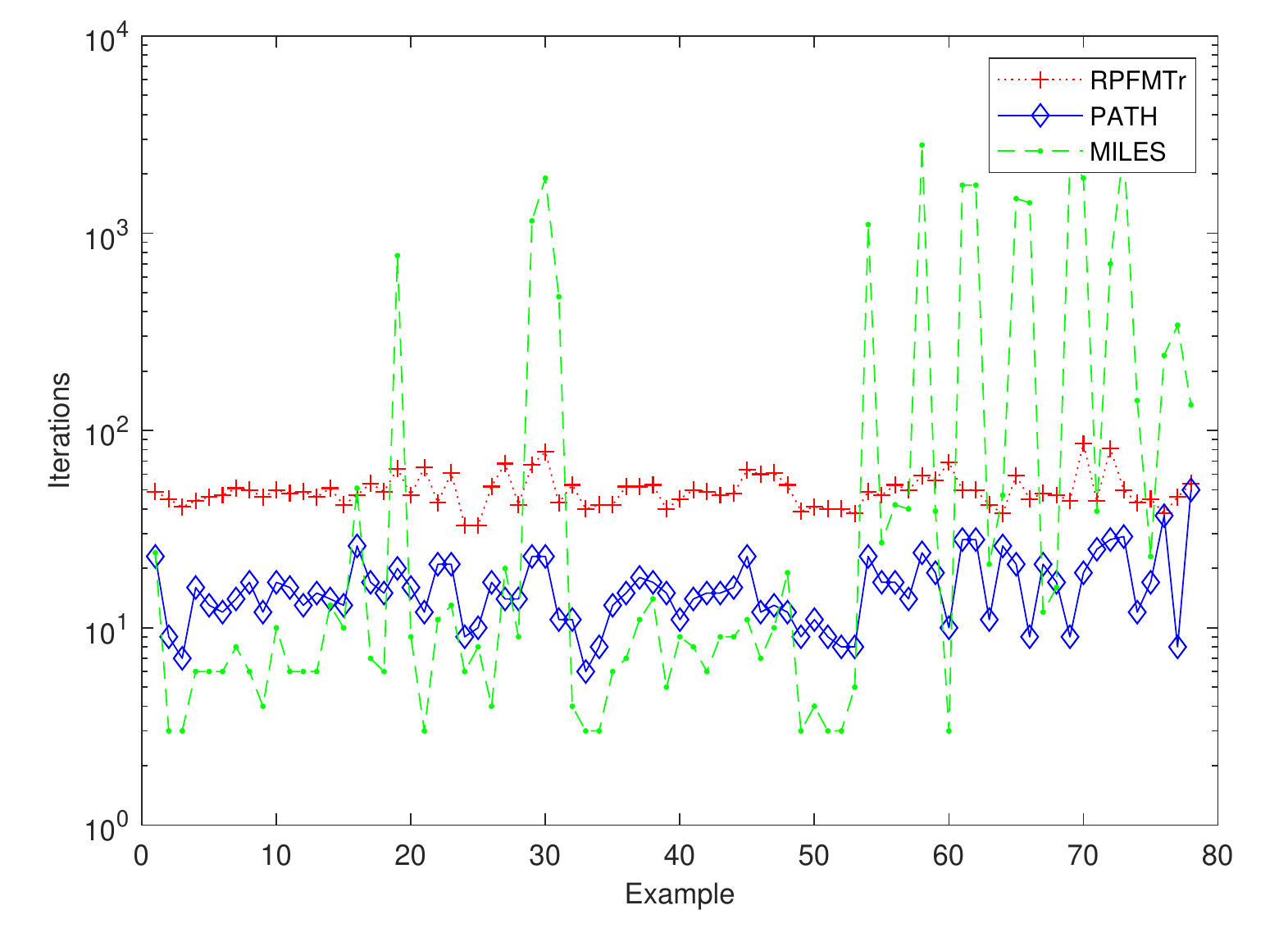}
        \caption{Iterations of RPFMTr, \,PATH and MILES for sparse LCPs.}
        \label{fig:ITSLCP}
\end{figure}

\begin{figure}[!htbp]
      \centering
        \includegraphics[width=0.80\textwidth]{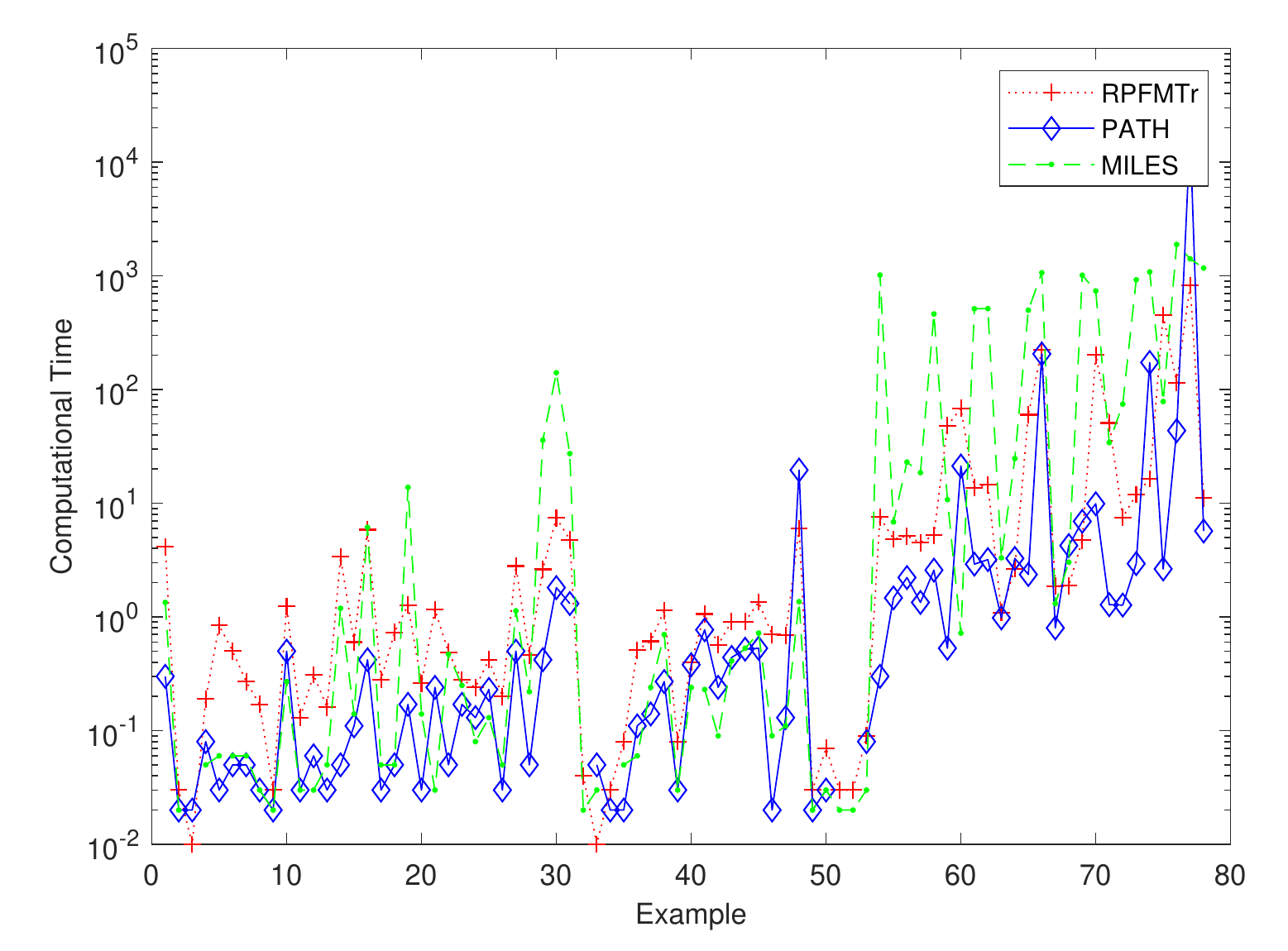}
        \caption{CPU time of RPFMTr, \,PATH and MILES for sparse LCPs.}
        \label{fig:CPUSLCP}
\end{figure}

\begin{figure}[!htbp]
      \centering
        \includegraphics[width=0.80\textwidth]{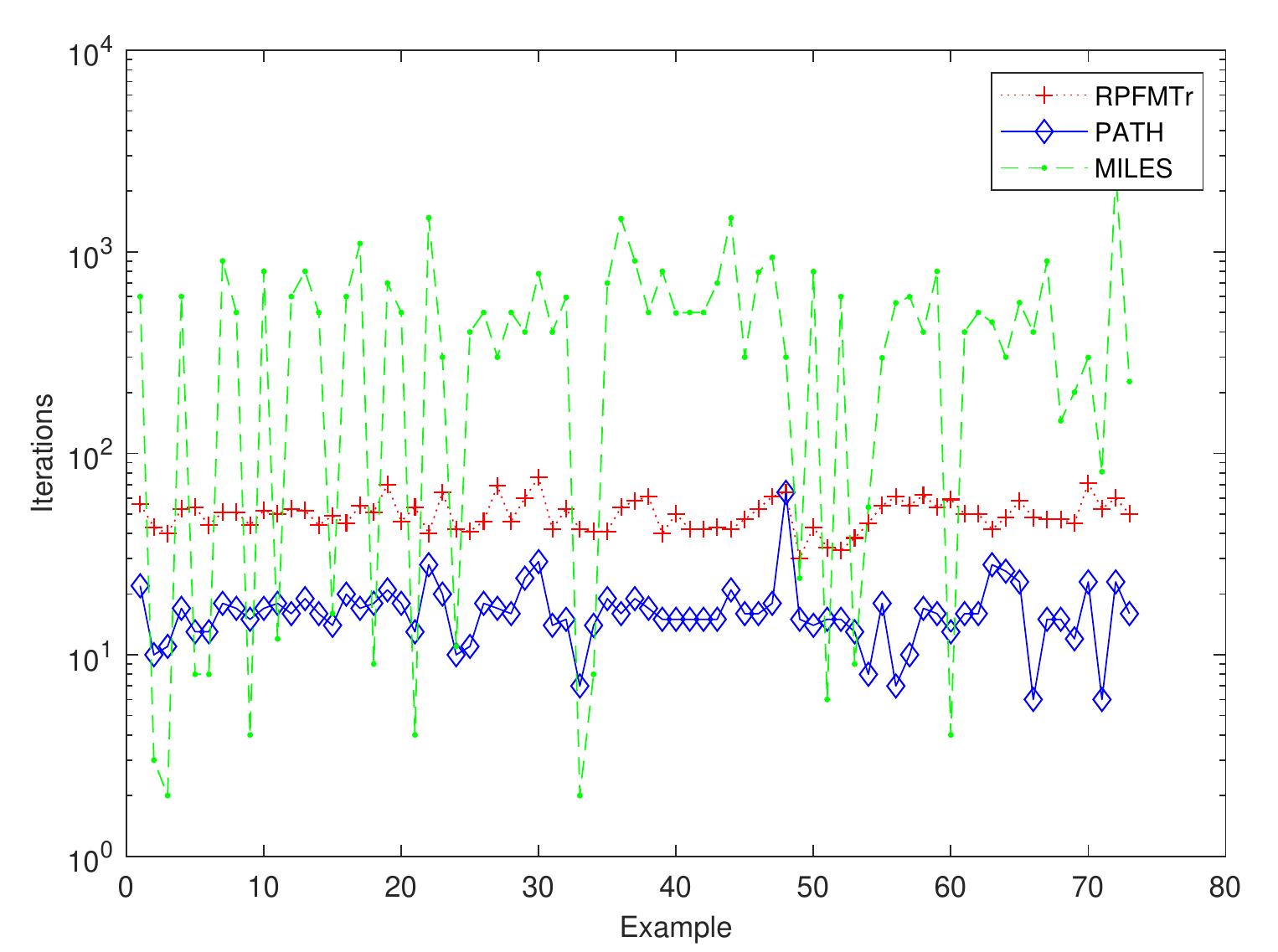}
        \caption{Iterations of RPFMTr, \,PATH and MILES for dense LCPs.}
        \label{fig:ITDLCP}
\end{figure}

\begin{figure}[!htbp]
      \centering
        \includegraphics[width=0.80\textwidth]{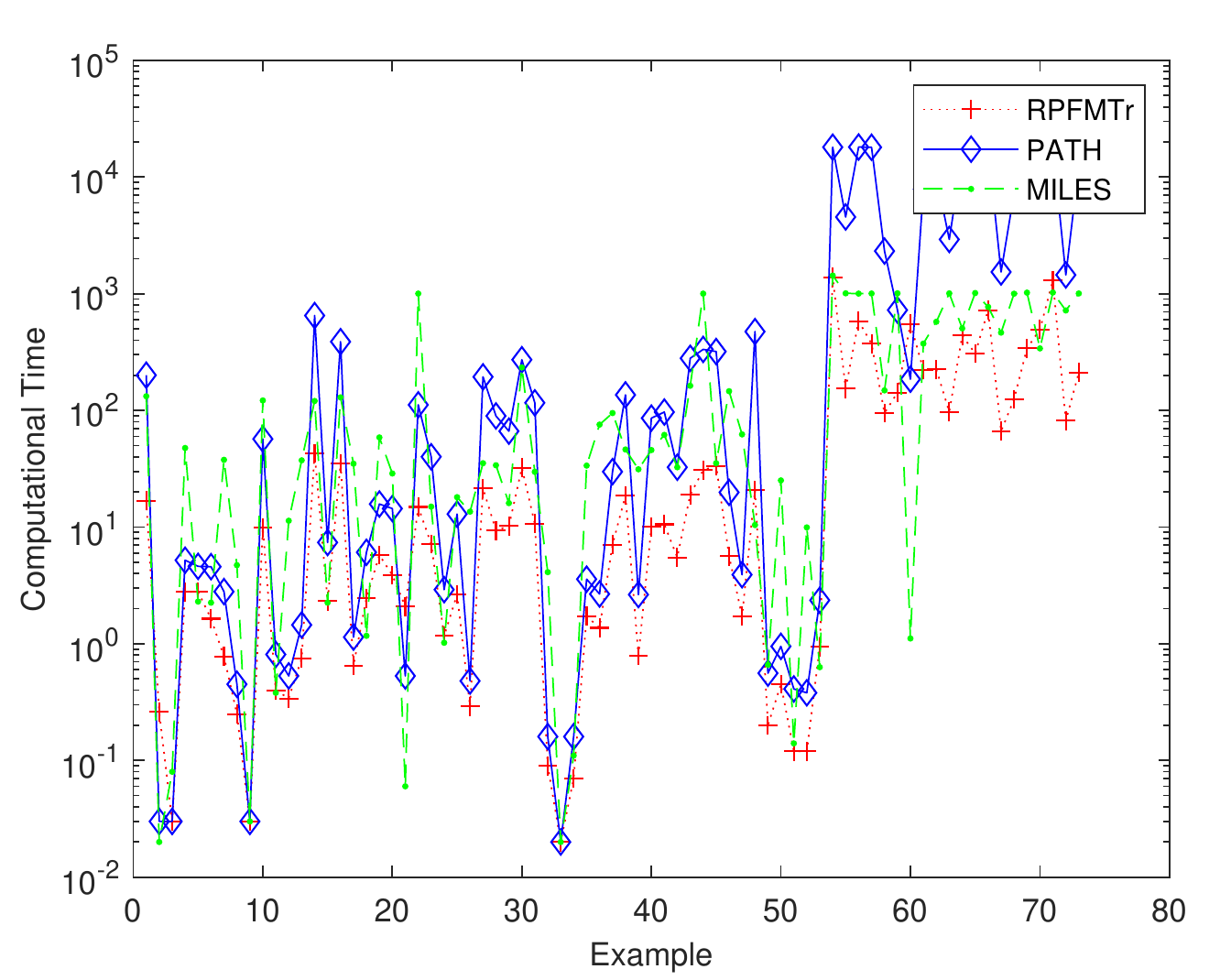}
        \caption{CPU time of RPFMTr, \,PATH and MILES for dense LCPs.}
        \label{fig:CPUDLCP}
\end{figure}

Numerical results of the sparse test problems are arranged in Tables
\ref{TABSLCP1}-\ref{TABSLCP3} and Figures \ref{fig:ITSLCP}-\ref{fig:CPUSLCP}.
And numerical results of the dense test problems are arranged in Tables
\ref{TABDLCP1}-\ref{TABDLCP3} and Figures \ref{fig:ITDLCP}-\ref{fig:CPUDLCP}.
``major'' in the fourth column of Tables \ref{TABSLCP1}-\ref{TABDLCP3} represents
the number of the linear mixed complementarity problems solved by a pivotal
Lemke's method of PATH \cite{FM2022}. ``minor'' in the fourth column of Tables
\ref{TABSLCP1}-\ref{TABDLCP3} represents the the number of pivots performed per
major iteration of PATH \cite{FM2022}. ``grad'' in the fourth column of Tables
\ref{TABSLCP1}-\ref{TABDLCP3} represents the cumulative number of
Jacobian evaluations used in PATH \cite{FM2022}.  ``major'' in the sixth column
of Tables \ref{TABSLCP1}-\ref{TABDLCP3} represents the number of Newton iterations
of MILES \cite{Rutherford2022}. ``pivots'' in the sixth column of Tables
\ref{TABSLCP1}-\ref{TABDLCP3} represents the number of Lemke pivots of MILES
\cite{Rutherford2022}. ``refactor'' in the sixth column of Tables
\ref{TABSLCP1}-\ref{TABDLCP3} represents the number of re-factorizations in
the LUSOL solver \cite{GMSW1991}, which is called by MILES for solving the linear
systems of equations.

\vskip 2mm

From Tables \ref{TABSLCP1}-\ref{TABDLCP3} and Figures
\ref{fig:ITSLCP}-\ref{fig:CPUDLCP}, we find that RPFMTr and PATH work well
for the sparse LCPs and the dense LCPs. However, MILES is not robust for
solving the sparse LCPs or the dense LCPs. For the sparse LCPs, PATH performs
better than RPFMTr and MILES from Tables \ref{TABSLCP1}-\ref{TABSLCP3} and
Figures \ref{fig:ITSLCP}-\ref{fig:CPUSLCP}. For the dense LCPs, from
Tables \ref{TABDLCP1}-\ref{TABDLCP3}, we find that PATH and MILES fail to
solve $5$ problems and $57$ problems of 73 test problems, respectively.
RPFMTr solves all the sparse test LCPs and the dense test LCPs.
Furthermore, from Tables \ref{TABDLCP1}-\ref{TABDLCP3} and Figures
\ref{fig:ITDLCP}-\ref{fig:CPUDLCP}, we find that the computational time of
RPFMTr is about $1/3$ to $1/10$ of that of PATH for the dense LCPs. Therefore,
RPFMTr is a robust and efficient solver for the LCPs, especially for the dense LCPs.

\newpage

\begin{table}[htbp]
  \newcommand{\tabincell}[2]{\begin{tabular}{@{}#1@{}}#2\end{tabular}}
  \scriptsize
  \centering
  \fontsize{8}{8}\selectfont
  \caption{Numerical results of small-scale sparse LCPs (no. 1-27).}
  \label{TABSLCP1}
  \resizebox{\textwidth}{!}{
  \begin{tabular}{|c|c|c|c|c|c|c|c|c|c|}
  \hline
  \multirow{2}{*}{Problems} & \multicolumn{2}{c|}{RPFMTr} & \multicolumn{2}{c|}{PATH} & \multicolumn{2}{c|}{MILES}  \\ \cline{2-7}
                \tabincell{c}{(n*n)}
                & \tabincell{c}{steps \\(time)} & \tabincell{c}{Terr\\}
                & \tabincell{c}{major+minor+grad \\(time)} & \tabincell{c}{Terr\\}
                & \tabincell{c}{major+pivots+refactor \\(time)} & \tabincell{c}{Terr\\} \\ \hline
  \tabincell{c}{Exam. 1 lp\_25fv47 \\ (n = 2697)}
  & \tabincell{c}{49 \\ (4.13)}  &\tabincell{c}{8.84e-07}
  & \tabincell{c}{8+949+23 \\ (0.30)} &\tabincell{c}{5.82e-08}
  & \tabincell{c}{1+1684+24 \\ (1.34)} &\tabincell{c}{7.48e-12} \\ \hline

  \tabincell{c}{Exam. 2 lp\_adlittle \\ (n = 194)}
  & \tabincell{c}{45 \\ (0.03)} &\tabincell{c}{9.61e-07}
  & \tabincell{c}{4+80+9 \\ (0.02)} &\tabincell{c}{4.44e-12 }
  & \tabincell{c}{1+1+3 \\ (0.02)} &\tabincell{c}{4.65e-14} \\ \hline

  \tabincell{c}{Exam. 3 lp\_afiro  \\ (n = 78)}
  & \tabincell{c}{41 \\ (0.01)} &\tabincell{c}{6.20e-07}
  & \tabincell{c}{3+5+7 \\ (0.02)} &\tabincell{c}{1.01e-10}
  & \tabincell{c}{1+1+3 \\ (0.00)} &\tabincell{c}{3.91e-14} \\ \hline

  \tabincell{c}{Exam. 4 lp\_agg \\ (n = 1103)}
  & \tabincell{c}{44 \\ (0.19)} &\tabincell{c}{9.85e-07}
  & \tabincell{c}{11+413+16 \\ (0.08)} &\tabincell{c}{2.87e-09}
  & \tabincell{c}{1+1+6 \\ (0.05)} &\tabincell{c}{3.74e-12} \\ \hline

  \tabincell{c}{Exam. 5 lp\_agg2 \\ (n = 1274)}
  & \tabincell{c}{46 \\ (0.84)} &\tabincell{c}{8.24e-07}
  & \tabincell{c}{10+303+13 \\ (0.03)} &\tabincell{c}{6.96e-10}
  & \tabincell{c}{1+959+6 \\ (0.06)} &\tabincell{c}{5.08e-12} \\ \hline

  \tabincell{c}{Exam. 6 lp\_agg3 \\ (n = 1274)}
  & \tabincell{c}{47 \\ (0.50)} &\tabincell{c}{6.64e-07}
  & \tabincell{c}{9+259+12 \\ (0.05)} &\tabincell{c}{2.41e-07 }
  & \tabincell{c}{1+955+6 \\ (0.06)} &\tabincell{c}{9.37e-12} \\ \hline

  \tabincell{c}{Exam. 7 lp\_bandm \\ (n = 777)}
  & \tabincell{c}{51 \\ (0.27)} &\tabincell{c}{8.58e-07}
  & \tabincell{c}{11+613+14 \\ (0.05)} &\tabincell{c}{3.31e-07 }
  & \tabincell{c}{1+1+8 \\ (0.06)} &\tabincell{c}{1.19e-12} \\ \hline

  \tabincell{c}{Exam. 8 lp\_beaconfd \\ (n = 468)}
  & \tabincell{c}{50 \\ (0.17)} &\tabincell{c}{6.06e-07}
  & \tabincell{c}{8+124+17 \\ (0.03)} &\tabincell{c}{3.76e-09 }
  & \tabincell{c}{1+1+6 \\ (0.03)} &\tabincell{c}{5.46e-12} \\ \hline

  \tabincell{c}{Exam. 9 lp\_blend \\ (n = 188)}
  & \tabincell{c}{46 \\ (0.03)} &\tabincell{c}{6.12e-07}
  & \tabincell{c}{9+143+12 \\ (0.02)} &\tabincell{c}{6.50e-08 }
  & \tabincell{c}{1+1+4 \\ (0.02)} &\tabincell{c}{1.96e-12} \\ \hline

  \tabincell{c}{Exam. 10 lp\_bnl1 \\ (n = 2229)}
  & \tabincell{c}{50 \\ (1.24)} &\tabincell{c}{8.60e-07}
  & \tabincell{c}{12+1727+17 \\ (0.50)} &\tabincell{c}{5.22e-11 }
  & \tabincell{c}{1+1+10 \\ (0.27)} &\tabincell{c}{1.62e-12} \\ \hline

  \tabincell{c}{Exam. 11 lp\_bore3d \\ (n = 567)}
  & \tabincell{c}{48 \\ (0.13)} &\tabincell{c}{8.63e-07}
  & \tabincell{c}{11+423+16 \\ (0.03)} &\tabincell{c}{6.99e-07 }
  & \tabincell{c}{1+1+6 \\ (0.03)} &\tabincell{c}{1.31e-12} \\ \hline

  \tabincell{c}{Exam. 12 lp\_brandy \\ (n = 523)}
  & \tabincell{c}{49 \\ (0.31)} &\tabincell{c}{7.42e-07}
  & \tabincell{c}{10+384+13 \\ (0.06)} &\tabincell{c}{2.84e-10 }
  & \tabincell{c}{1+1+6 \\ (0.03)} &\tabincell{c}{1.31e-12} \\ \hline

  \tabincell{c}{Exam. 13 lp\_capri \\ (n = 753)}
  & \tabincell{c}{46 \\ (0.16)} &\tabincell{c}{9.07e-07}
  & \tabincell{c}{11+573+15 \\ (0.03)} &\tabincell{c}{7.18e-07}
  & \tabincell{c}{1+1+6 \\ (0.05)} &\tabincell{c}{2.60e-12} \\ \hline

  \tabincell{c}{Exam. 14 lp\_czprob \\ (n = 4491)}
  & \tabincell{c}{51 \\ (3.37)} &\tabincell{c}{8.74e-07}
  & \tabincell{c}{4+25+14 \\ (0.05)} &\tabincell{c}{1.90e-10 }
  & \tabincell{c}{1+1+13 \\ (1.19)} &\tabincell{c}{1.82e-12} \\ \hline

  \tabincell{c}{Exam. 15 lp\_degen2 \\ (n = 1201)}
  & \tabincell{c}{42 \\ (0.60)} &\tabincell{c}{8.53e-07}
  & \tabincell{c}{4+447+13 \\ (0.11)} &\tabincell{c}{1.71e-07 }
  & \tabincell{c}{1+1+10 \\ (0.14)} &\tabincell{c}{2.22e-16} \\ \hline

  \tabincell{c}{Exam. 16 lp\_degen3 \\ (n = 4107)}
  & \tabincell{c}{47 \\ (5.85)} &\tabincell{c}{6.85e-07}
  & \tabincell{c}{4+185+26 \\ (0.42)} &\tabincell{c}{3.00e-07 }
  & \tabincell{c}{1+1+51 \\ (6.13)} &\tabincell{c}{0.00e-00} \\ \hline

  \tabincell{c}{Exam. 17 lp\_e226 \\ (n = 695)}
  & \tabincell{c}{54 \\ (0.28)} &\tabincell{c}{9.05e-07}
  & \tabincell{c}{9+412+17 \\ (0.03)} &\tabincell{c}{2.94e-07 }
  & \tabincell{c}{1+1+7 \\ (0.05)} &\tabincell{c}{5.86e-12} \\ \hline

  \tabincell{c}{Exam. 18 lp\_etamacro \\ (n = 1216)}
  & \tabincell{c}{49 \\ (0.73)} &\tabincell{c}{6.31e-07}
  & \tabincell{c}{6+440+15 \\ (0.05)} &\tabincell{c}{1.87e-12 }
  & \tabincell{c}{1+58+6 \\ (0.05)} &\tabincell{c}{6.82e-12} \\ \hline

  \tabincell{c}{Exam. 19 lp\_fffff800 \\ (n = 1552)}
  & \tabincell{c}{64 \\ (1.27)} &\tabincell{c}{6.05e-07}
  & \tabincell{c}{15+1111+20 \\ (0.17)} &\tabincell{c}{1.39e-07 }
  &\tabincell{c}{\textcolor{red}{100+72624+770} \\ \textcolor{red}{(13.80)}} &\tabincell{c}{\textcolor{red}{1.78e+01}  \\ \textcolor{red}{(failed)}} \\ \hline

  \tabincell{c}{Exam. 20 lp\_finnis \\ (n = 1561)}
  & \tabincell{c}{47 \\ (0.26)} &\tabincell{c}{6.06e-07}
  & \tabincell{c}{8+313+16 \\ (0.03)} &\tabincell{c}{5.24e-09 }
  & \tabincell{c}{1+1+9 \\ (0.14)} &\tabincell{c}{4.55e-13} \\ \hline

  \tabincell{c}{Exam. 21 lp\_fit1d \\ (n = 1073)}
  & \tabincell{c}{65 \\ (1.16)} &\tabincell{c}{6.45e-07}
  & \tabincell{c}{8+671+12 \\ (0.24)} &\tabincell{c}{2.34e-11 }
  & \tabincell{c}{1+1+3 \\ (0.03)} &\tabincell{c}{1.75e-10} \\ \hline

  \tabincell{c}{Exam. 22 lp\_ganges \\ (n = 3015)}
  & \tabincell{c}{43 \\ (0.49)} &\tabincell{c}{8.88e-07}
  & \tabincell{c}{4+98+21 \\ (0.05)} &\tabincell{c}{4.95e-07 }
  & \tabincell{c}{1+1+11 \\ (0.47)} &\tabincell{c}{2.40e-14} \\ \hline

  \tabincell{c}{Exam. 23 lp\_gfrd\_pnc \\ (n = 1776)}
  & \tabincell{c}{61 \\ (0.28)} &\tabincell{c}{7.85e-07}
  & \tabincell{c}{18+1250+21 \\ (0.17)} &\tabincell{c}{3.80e-11 }
  & \tabincell{c}{1+716+13 \\ (0.25)} &\tabincell{c}{3.64e-12} \\ \hline

  \tabincell{c}{Exam. 24 lp\_grow15 \\ (n = 945)}
  & \tabincell{c}{33 \\ (0.24)} &\tabincell{c}{7.23e-07}
  & \tabincell{c}{7+703+9 \\ (0.13)} &\tabincell{c}{8.66e-07 }
  & \tabincell{c}{1+927+6 \\ (0.08)} &\tabincell{c}{3.33e-14} \\ \hline

  \tabincell{c}{Exam. 25 lp\_grow22 \\ (n = 1386)}
  & \tabincell{c}{33 \\ (0.42)} &\tabincell{c}{7.27e-07}
  & \tabincell{c}{8+1054+10 \\ (0.23)} &\tabincell{c}{3.86e-11}
  & \tabincell{c}{1+1403+8 \\ (0.13)} &\tabincell{c}{1.20e-14} \\ \hline

  \tabincell{c}{Exam. 26 lp\_lotfi \\ (n = 519)}
  & \tabincell{c}{52 \\ (0.20)} &\tabincell{c}{6.05e-07}
  & \tabincell{c}{13+263+17 \\ (0.03)} &\tabincell{c}{8.42e-07 }
  & \tabincell{c}{1+1+4 \\ (0.05)} &\tabincell{c}{9.10e-13} \\ \hline

  \tabincell{c}{Exam. 27 lp\_maros \\ (n = 2812)}
  & \tabincell{c}{68 \\ (2.80)} &\tabincell{c}{6.18e-07}
  & \tabincell{c}{13+3222+14 \\ (0.50)} &\tabincell{c}{2.16e-09 }
  & \tabincell{c}{1+776+20 \\ (1.13)} &\tabincell{c}{2.95e-09} \\ \hline

\end{tabular}}%
\end{table}%

\begin{table}[htbp]
  \newcommand{\tabincell}[2]{\begin{tabular}{@{}#1@{}}#2\end{tabular}}
  \scriptsize
  \centering
  \fontsize{8}{8}\selectfont
  \caption{Numerical results of small-scale sparse LCPs (no. 28-53).}
  \label{TABSLCP2}
  \resizebox{\textwidth}{!}{
  \begin{tabular}{|c|c|c|c|c|c|c|c|c|c|}
  \hline
  \multirow{2}{*}{Problems} & \multicolumn{2}{c|}{RPFMTr} & \multicolumn{2}{c|}{PATH} & \multicolumn{2}{c|}{MILES}  \\ \cline{2-7}
              \tabincell{c}{(n*n)}
                & \tabincell{c}{steps \\(time)} & \tabincell{c}{Terr\\}
                & \tabincell{c}{major+minor+grad \\(time)} & \tabincell{c}{Terr\\}
                & \tabincell{c}{major+pivots+refactor \\(time)} & \tabincell{c}{Terr\\} \\ \hline
  \tabincell{c}{Exam. 28 lp\_modszk1 \\ (n = 2307)}
  & \tabincell{c}{42 \\ (0.46)} &\tabincell{c}{8.61e-07}
  & \tabincell{c}{5+370+14 \\ (0.05)} &\tabincell{c}{7.91e-11 }
  & \tabincell{c}{1+1+9 \\ (0.22)} &\tabincell{c}{8.70e-14} \\ \hline

  \tabincell{c}{Exam. 29 lp\_perold \\ (n = 2131)}
  & \tabincell{c}{67 \\ (2.62)} &\tabincell{c}{9.15e-07}
  & \tabincell{c}{12+1498+23 \\ (0.42)} &\tabincell{c}{1.72e-08 }
  &\tabincell{c}{\textcolor{red}{100+18168+1156} \\ \textcolor{red}{(35.88)}} &\tabincell{c}{\textcolor{red}{6.64e+01}  \\ \textcolor{red}{(failed)}} \\ \hline

  \tabincell{c}{Exam. 30 lp\_pilot\_ja \\ (n = 3207)}
  & \tabincell{c}{78 \\ (7.43)} &\tabincell{c}{7.94e-07}
  & \tabincell{c}{14+2468+23 \\ (1.81)} &\tabincell{c}{3.29e-07 }
  &\tabincell{c}{\textcolor{red}{100+171048+1900} \\ \textcolor{red}{(140.47)}} &\tabincell{c}{\textcolor{red}{3.92e+04}  \\ \textcolor{red}{(failed)}} \\ \hline

  \tabincell{c}{Exam. 31 lp\_qap8 \\ (n = 2544)}
  & \tabincell{c}{43 \\ (4.78)} &\tabincell{c}{8.80e-07}
  & \tabincell{c}{9+2413+11 \\ (1.31)} &\tabincell{c}{1.83e-11 }
  & \tabincell{c}{6+50530+477 \\ (27.30)} &\tabincell{c}{2.13e-14 } \\ \hline

  \tabincell{c}{Exam. 32 lp\_lp\_recipe \\ (n = 295)}
  & \tabincell{c}{53 \\ (0.04)} &\tabincell{c}{7.81e-07}
  & \tabincell{c}{8+103+11 \\ (0.00)} &\tabincell{c}{7.33e-12 }
  & \tabincell{c}{1+1+4 \\ (0.02)} &\tabincell{c}{5.12e-13} \\ \hline

  \tabincell{c}{Exam. 33 lp\_sc50a \\ (n = 128)}
  & \tabincell{c}{40 \\ (0.01)} &\tabincell{c}{8.30e-07}
  & \tabincell{c}{3+59+6 \\ (0.05)} &\tabincell{c}{9.21e-07 }
  & \tabincell{c}{1+1+3 \\ (0.03)} &\tabincell{c}{1.33e-15} \\ \hline

  \tabincell{c}{Exam. 34 lp\_scagr7 \\ (n = 314)}
  & \tabincell{c}{42 \\ (0.03)} &\tabincell{c}{6.28e-07}
  & \tabincell{c}{5+110+8 \\ (0.02)} &\tabincell{c}{1.41e-10 }
  & \tabincell{c}{1+1+3 \\ (0.00)} &\tabincell{c}{1.07e-14} \\ \hline

  \tabincell{c}{Exam. 35 lp\_scagr25 \\ (n = 1142)}
  & \tabincell{c}{42 \\ (0.08)} &\tabincell{c}{6.30e-07}
  & \tabincell{c}{4+205+13 \\ (0.02)} &\tabincell{c}{4.24e-07 }
  & \tabincell{c}{1+1+6 \\ (0.05)} &\tabincell{c}{1.07e-14} \\ \hline

  \tabincell{c}{Exam. 36 lp\_scfxm1 \\ (n = 930)}
  & \tabincell{c}{52 \\ (0.51)} &\tabincell{c}{9.05e-07}
  & \tabincell{c}{12+721+15 \\ (0.11)} &\tabincell{c}{7.23e-10 }
  & \tabincell{c}{1+1+7 \\ (0.06)} &\tabincell{c}{2.62e-12} \\ \hline

  \tabincell{c}{Exam. 37 lp\_scfxm2 \\ (n = 1860)}
  & \tabincell{c}{52 \\ (0.61)} &\tabincell{c}{8.84e-07}
  & \tabincell{c}{14+1357+18 \\ (0.14)} &\tabincell{c}{3.36e-09 }
  & \tabincell{c}{1+1+11 \\ (0.24)} &\tabincell{c}{1.00e-11} \\ \hline

  \tabincell{c}{Exam. 38 lp\_scfxm3 \\ (n = 2790)}
  & \tabincell{c}{53 \\ (1.14)} &\tabincell{c}{6.17e-07}
  & \tabincell{c}{14+2110+17 \\ (0.27)} &\tabincell{c}{4.96e-07 }
  & \tabincell{c}{1+200+14 \\ (0.70)} &\tabincell{c}{7.38e-12} \\ \hline

  \tabincell{c}{Exam. 39 lp\_scorpion \\ (n = 854)}
  & \tabincell{c}{40 \\ (0.08)} &\tabincell{c}{8.40e-07}
  & \tabincell{c}{5+111+15 \\ (0.03)} &\tabincell{c}{1.35e-07 }
  & \tabincell{c}{1+1+5 \\ (0.03)} &\tabincell{c}{6.22e-15} \\ \hline

  \tabincell{c}{Exam. 40 lp\_shell \\ (n = 2313)}
  & \tabincell{c}{45 \\ (0.40)} &\tabincell{c}{6.34e-07}
  & \tabincell{c}{9+2448+11 \\ (0.38)} &\tabincell{c}{5.86e-10 }
  & \tabincell{c}{1+1+9 \\ (0.24)} &\tabincell{c}{0.00e-00} \\ \hline

  \tabincell{c}{Exam. 41 lp\_ship04l \\ (n = 2568)}
  & \tabincell{c}{50 \\ (1.06)} &\tabincell{c}{7.03e-07}
  & \tabincell{c}{11+3359+14 \\ (0.77)} &\tabincell{c}{4.05e-07 }
  & \tabincell{c}{1+1+8\\ (0.23)} &\tabincell{c}{1.81e-14} \\ \hline

  \tabincell{c}{Exam. 42 lp\_ship04s \\ (n = 1908)}
  & \tabincell{c}{49 \\ (0.57)} &\tabincell{c}{9.60e-07}
  & \tabincell{c}{12+2376+15 \\ (0.24)} &\tabincell{c}{4.04e-07 }
  & \tabincell{c}{1+1+6 \\ (0.09)} &\tabincell{c}{1.81e-14} \\ \hline

  \tabincell{c}{Exam. 43 lp\_ship08s \\ (n = 3245)}
  & \tabincell{c}{47 \\ (0.91)} &\tabincell{c}{5.31e-07}
  & \tabincell{c}{12+3328+15 \\ (0.44)} &\tabincell{c}{9.32e-07 }
  & \tabincell{c}{1+1+9 \\ (0.41)} &\tabincell{c}{5.68e-14} \\ \hline

  \tabincell{c}{Exam. 44 lp\_ship12s \\ (n = 4020)}
  & \tabincell{c}{48 \\ (0.91)} &\tabincell{c}{2.68e-07}
  & \tabincell{c}{13+3446+16 \\ (0.52)} &\tabincell{c}{3.14e-07 }
  & \tabincell{c}{1+1+9 \\ (0.53)} &\tabincell{c}{4.17e-14} \\ \hline

  \tabincell{c}{Exam. 45 lp\_sierra \\ (n = 3962)}
  & \tabincell{c}{63 \\ (1.35)} &\tabincell{c}{6.79e-07}
  & \tabincell{c}{12+2135+23 \\ (0.53)} &\tabincell{c}{7.90e-07 }
  & \tabincell{c}{1+1+11 \\ (0.72)} &\tabincell{c}{4.37e-11} \\ \hline

  \tabincell{c}{Exam. 46 lp\_standgub \\ (n = 1744)}
  & \tabincell{c}{60 \\ (0.70)} &\tabincell{c}{9.29e-07}
  & \tabincell{c}{4+99+12 \\ (0.02)} &\tabincell{c}{2.88e-12 }
  & \tabincell{c}{1+1+7 \\ (0.09)} &\tabincell{c}{1.16e-13} \\ \hline

  \tabincell{c}{Exam. 47 lp\_tuff \\ (n = 961)}
  & \tabincell{c}{61 \\ (0.69)} &\tabincell{c}{3.85e-07}
  & \tabincell{c}{9+828+13 \\ (0.13)} &\tabincell{c}{4.81e-07 }
  & \tabincell{c}{1+1+10 \\ (0.11)} &\tabincell{c}{3.82e-11} \\ \hline

  \tabincell{c}{Exam. 48 lp\_wood1p \\ (n = 2839)}
  & \tabincell{c}{53 \\ (6.00)} &\tabincell{c}{7.48e-07}
  & \tabincell{c}{7+1285+12 \\ (19.56)} &\tabincell{c}{9.65e-10 }
  & \tabincell{c}{1+3237+19 \\ (1.36)} &\tabincell{c}{5.09e-09} \\ \hline

  \tabincell{c}{Exam. 49 lpi\_box1 \\ (n = 492)}
  & \tabincell{c}{39 \\ (0.03)} &\tabincell{c}{6.17e-07}
  & \tabincell{c}{7+290+9 \\ (0.02)} &\tabincell{c}{7.28e-07 }
  & \tabincell{c}{1+279+3 \\ (0.02)} &\tabincell{c}{1.11e-16} \\ \hline

  \tabincell{c}{Exam. 50 lpi\_cplex2 \\ (n = 602)}
  & \tabincell{c}{41 \\ (0.07)} &\tabincell{c}{7.33e-07}
  & \tabincell{c}{9+422+11 \\ (0.03)} &\tabincell{c}{1.45e-09 }
  & \tabincell{c}{1+1+4 \\ (0.03)} &\tabincell{c}{4.44e-15} \\ \hline

  \tabincell{c}{Exam. 51 lpi\_ex72a \\ (n = 412)}
  & \tabincell{c}{40 \\ (0.03)} &\tabincell{c}{6.16e-07}
  & \tabincell{c}{7+164+9 \\ (0.00)} &\tabincell{c}{1.24e-12 }
  & \tabincell{c}{1+1+3 \\ (0.02)} &\tabincell{c}{0.00e-00} \\ \hline

  \tabincell{c}{Exam. 52 lpi\_ex73a \\ (n = 404)}
  & \tabincell{c}{40 \\ (0.03)} &\tabincell{c}{7.56e-07}
  & \tabincell{c}{6+125+8 \\ (0.00)} &\tabincell{c}{6.75e-10 }
  & \tabincell{c}{1+1+3 \\ (0.02)} &\tabincell{c}{0.00e-00} \\ \hline

  \tabincell{c}{Exam. 53 lpi\_mondou2 \\ (n = 916)}
  & \tabincell{c}{38 \\ (0.09)} &\tabincell{c}{6.85e-07}
  & \tabincell{c}{6+555+8 \\ (0.08)} &\tabincell{c}{5.02e-07 }
  & \tabincell{c}{1+1+5 \\ (0.03)} &\tabincell{c}{0.00e-00} \\ \hline

\end{tabular}}%
\end{table}%

\begin{table}[htbp]
  \newcommand{\tabincell}[2]{\begin{tabular}{@{}#1@{}}#2\end{tabular}}
  \scriptsize
  \centering
  \fontsize{8}{8}\selectfont
  \caption{Numerical results of large-scale sparse LCPs (no. 54-78).}
  \label{TABSLCP3}
  \resizebox{\textwidth}{!}{
  \begin{tabular}{|c|c|c|c|c|c|c|c|c|c|}
  \hline
  \multirow{2}{*}{Problems} & \multicolumn{2}{c|}{RPFMTr} & \multicolumn{2}{c|}{PATH} & \multicolumn{2}{c|}{MILES}  \\ \cline{2-7}
              \tabincell{c}{(n*n)}
                & \tabincell{c}{steps \\(time)} & \tabincell{c}{Terr\\}
                & \tabincell{c}{major+minor+grad \\(time)} & \tabincell{c}{Terr\\}
                & \tabincell{c}{major+pivots+refactor \\(time)} & \tabincell{c}{Terr\\} \\ \hline

  \tabincell{c}{Exam. 54 lp\_80bau3b \\ (n = 14323)}
  & \tabincell{c}{49 \\ (7.56)} &\tabincell{c}{6.82e-07}
  & \tabincell{c}{9+858+23 \\ (0.30)} &\tabincell{c}{1.52e-10}
  &\tabincell{c}{\textcolor{red}{41+101114+1107} \\ \textcolor{red}{(1015.16)}} &\tabincell{c}{\textcolor{red}{6.19e+01}  \\ \textcolor{red}{(failed)}}  \\ \hline

  \tabincell{c}{Exam. 55 lp\_bnl2 \\ (n = 6810)}
  & \tabincell{c}{47 \\ (4.83)} &\tabincell{c}{7.01e-07}
  & \tabincell{c}{12+4597+17 \\ (1.47)} &\tabincell{c}{2.72e-11 }
  & \tabincell{c}{1+1+27 \\ (6.84)} &\tabincell{c}{1.24e-12} \\ \hline

  \tabincell{c}{Exam. 56 lp\_cre\_a  \\ (n = 10764)}
  & \tabincell{c}{53 \\ (5.14)} &\tabincell{c}{8.39e-07}
  & \tabincell{c}{12+5755+17 \\ (2.22)} &\tabincell{c}{4.09e-07 }
  & \tabincell{c}{1+1+42 \\ (22.98)} &\tabincell{c}{1.48e-12} \\ \hline

  \tabincell{c}{Exam. 57 lp\_cre\_c \\ (n = 9479)}
  & \tabincell{c}{50 \\ (4.50)} &\tabincell{c}{8.58e-07}
  & \tabincell{c}{11+3604+14 \\ (1.34)} &\tabincell{c}{1.77e-09 }
  & \tabincell{c}{1+490+40 \\ (18.56)} &\tabincell{c}{1.59e-12} \\ \hline

  \tabincell{c}{Exam. 58 lp\_cycle\\ (n = 5274)}
  & \tabincell{c}{59 \\ (5.21)} &\tabincell{c}{4.59e-07}
  & \tabincell{c}{17+4532+24 \\ (2.58)} &\tabincell{c}{8.67e-09 }
  &\tabincell{c}{\textcolor{red}{100+197608+2800} \\ \textcolor{red}{(461.45)}} &\tabincell{c}{\textcolor{red}{7.15e+03}  \\ \textcolor{red}{(failed)}} \\ \hline

  \tabincell{c}{Exam. 59 lp\_d6cube \\ (n = 6599)}
  & \tabincell{c}{56 \\ (48.11)} &\tabincell{c}{2.56e-07}
  & \tabincell{c}{5+941+19 \\ (0.53)} &\tabincell{c}{7.37e-07 }
  & \tabincell{c}{1+6085+39 \\ (10.75)} &\tabincell{c}{2.73e-12} \\ \hline

  \tabincell{c}{Exam. 60 lp\_fit2d \\ (n = 10549)}
  & \tabincell{c}{69 \\ (68.27)} &\tabincell{c}{6.72e-07}
  & \tabincell{c}{7+2114+10 \\ (21.23)} &\tabincell{c}{1.53e-09 }
  & \tabincell{c}{1+1+3 \\ (0.72)} &\tabincell{c}{3.49e-10} \\ \hline

  \tabincell{c}{Exam. 61 lp\_greenbea \\ (n = 7990)}
  & \tabincell{c}{50 \\ (13.68)} &\tabincell{c}{7.29e-07}
  & \tabincell{c}{12+3899+28 \\ (2.91)} &\tabincell{c}{1.01e-09 }
  &\tabincell{c}{\textcolor{red}{100+74070+1754} \\ \textcolor{red}{(513.08)}} &\tabincell{c}{\textcolor{red}{9.90e+01}  \\ \textcolor{red}{(failed)}}  \\ \hline

  \tabincell{c}{Exam. 62 lp\_greenbeb \\ (n = 7990)}
  & \tabincell{c}{50 \\ (14.48)} &\tabincell{c}{7.29e-07}
  & \tabincell{c}{12+3899+28 \\ (3.19)} &\tabincell{c}{1.01e-09 }
  &\tabincell{c}{\textcolor{red}{100+74070+1754} \\ \textcolor{red}{(514.13)}} &\tabincell{c}{\textcolor{red}{9.90e+01}  \\ \textcolor{red}{(failed)}} \\ \hline

  \tabincell{c}{Exam. 63 lp\_ken\_07 \\ (n = 6028)}
  & \tabincell{c}{42 \\ (1.09)} &\tabincell{c}{6.24e-07}
  & \tabincell{c}{8+4004+11 \\ (0.98)} &\tabincell{c}{9.71e-09 }
  & \tabincell{c}{1+1+21 \\ (3.31)} &\tabincell{c}{0.00e-00} \\ \hline

  \tabincell{c}{Exam. 64 lp\_pds\_02 \\ (n = 10669)}
  & \tabincell{c}{38 \\ (2.66)} &\tabincell{c}{9.65e-07}
  & \tabincell{c}{5+3264+26 \\ (3.27)} &\tabincell{c}{7.09e-09 }
  & \tabincell{c}{1+1+47 \\ (24.70)} &\tabincell{c}{0.00e-00} \\ \hline

  \tabincell{c}{Exam. 65 lp\_pilot87 \\ (n = 8710)}
  & \tabincell{c}{59 \\ (59.97)} &\tabincell{c}{7.96e-07}
  & \tabincell{c}{10+3002+21 \\ (2.36)} &\tabincell{c}{6.47e-07 }
  &\tabincell{c}{\textcolor{red}{100+121000+1500} \\ \textcolor{red}{(498.08)}} &\tabincell{c}{\textcolor{red}{1.00e+03}  \\ \textcolor{red}{(failed)}} \\ \hline

  \tabincell{c}{Exam. 66 lp\_qap12 \\ (n = 12048)}
  & \tabincell{c}{45 \\ (223.00)} &\tabincell{c}{6.25e-07}
  & \tabincell{c}{7+9138+9 \\ (205.47)} &\tabincell{c}{8.54e-07 }
  &\tabincell{c}{\textcolor{red}{12+145940+1427} \\ \textcolor{red}{(1064.42)}} &\tabincell{c}{\textcolor{red}{1.07e+01}  \\ \textcolor{red}{(failed)}} \\ \hline

  \tabincell{c}{Exam. 67 lp\_ship08l \\ (n = 5141)}
  & \tabincell{c}{48 \\ (1.87)} &\tabincell{c}{5.86e-07}
  & \tabincell{c}{10+1676+21 \\ (0.80)} &\tabincell{c}{5.38e-07 }
  & \tabincell{c}{1+1+12 \\ (1.31)} &\tabincell{c}{5.68e-14} \\ \hline

  \tabincell{c}{Exam. 68 lp\_ship12l \\ (n = 6684)}
  & \tabincell{c}{47 \\ (1.88)} &\tabincell{c}{8.51e-07}
  & \tabincell{c}{14+8950+17 \\ (4.24)} &\tabincell{c}{1.22e-10 }
  & \tabincell{c}{1+1+16 \\ (3.03)} &\tabincell{c}{4.17e-14} \\ \hline

  \tabincell{c}{Exam. 69 lp\_truss \\ (n = 9806)}
  & \tabincell{c}{44 \\ (4.72)} &\tabincell{c}{6.22e-07}
  & \tabincell{c}{6+5607+9 \\ (6.92)} &\tabincell{c}{4.99e-08 }
  &\tabincell{c}{\textcolor{red}{98+299152+2379} \\ \textcolor{red}{(1011.18)}} &\tabincell{c}{\textcolor{red}{1.67e+00}  \\ \textcolor{red}{(failed)}} \\ \hline

  \tabincell{c}{Exam. 70 lp\_woodw \\ (n = 9516)}
  & \tabincell{c}{86 \\ (201.33)} &\tabincell{c}{6.31e-07}
  & \tabincell{c}{10+4715+19 \\ (9.89)} &\tabincell{c}{1.78e-09 }
  &\tabincell{c}{\textcolor{red}{100+232954+1904} \\ \textcolor{red}{(735.08)}} &\tabincell{c}{\textcolor{red}{2.00e+05}  \\ \textcolor{red}{(failed)}} \\ \hline

  \tabincell{c}{Exam. 71 lpi\_bgindy \\ (n = 13551)}
  & \tabincell{c}{44 \\ (50.77)} &\tabincell{c}{8.49e-07}
  & \tabincell{c}{4+614+25 \\ (1.28)} &\tabincell{c}{1.20e-07 }
  & \tabincell{c}{1+1+39 \\ (34.14)} &\tabincell{c}{9.06e-08} \\ \hline

  \tabincell{c}{Exam. 72 lpi\_gran \\ (n = 5183)}
  & \tabincell{c}{81 \\ (7.50)} &\tabincell{c}{9.94e-07}
  & \tabincell{c}{15+1988+28 \\ (1.27)} &\tabincell{c}{8.14e-07 }
  &\tabincell{c}{\textcolor{red}{100+70036+700} \\ \textcolor{red}{(74.33)}} &\tabincell{c}{\textcolor{red}{5.68e+03}  \\ \textcolor{red}{(failed)}}  \\ \hline

  \tabincell{c}{Exam. 73 lpi\_greenbea \\ (n = 7989)}
  & \tabincell{c}{50 \\ (11.91)} &\tabincell{c}{9.94e-07}
  & \tabincell{c}{11+3173+29 \\ (2.94)} &\tabincell{c}{4.68e-07 }
  &\tabincell{c}{\textcolor{red}{100+303396+2731} \\ \textcolor{red}{(921.31)}} &\tabincell{c}{\textcolor{red}{6.72e+01}  \\ \textcolor{red}{(failed)}}  \\ \hline

  \tabincell{c}{Exam. 74 lp\_ken\_11 \\ (n = 36043)}
  & \tabincell{c}{43 \\ (16.37)} &\tabincell{c}{6.95e-07}
  & \tabincell{c}{9+23833+12 \\ (172.78)} &\tabincell{c}{2.65e-08 }
  & \tabincell{c}{1+1+142 \\ (1081.75)} &\tabincell{c}{0.00e-00} \\ \hline

  \tabincell{c}{Exam. 75 lp\_osa\_07 \\ (n = 26185)}
  & \tabincell{c}{45 \\ (451.93)} &\tabincell{c}{6.55e-07}
  & \tabincell{c}{4+954+17 \\ (2.64)} &\tabincell{c}{1.77e-11 }
  & \tabincell{c}{1+1+23 \\ (78.39)} &\tabincell{c}{3.87e-09} \\ \hline

  \tabincell{c}{Exam. 76 lp\_pds\_06 \\ (n = 39232)}
  & \tabincell{c}{38 \\ (115.16)} &\tabincell{c}{9.02e-07}
  & \tabincell{c}{5+6809+37 \\ (43.47)} &\tabincell{c}{9.00e-11 }
  & \tabincell{c}{1+1+240 \\ (1885.67)} &\tabincell{c}{0.00e-00} \\ \hline

  \tabincell{c}{Exam. 77 lp\_qap15 \\ (n = 28605)}
  & \tabincell{c}{46 \\ (824.18)} &\tabincell{c}{6.34e-07}
  &\tabincell{c}{\textcolor{red}{5+20058+8} \\ \textcolor{red}{(18004.75)}} &\tabincell{c}{\textcolor{red}{2.11e+01}  \\ \textcolor{red}{(failed)}}
  &\tabincell{c}{\textcolor{red}{3+39920+342} \\ \textcolor{red}{(1410.55)}} &\tabincell{c}{\textcolor{red}{1.20e+01}  \\ \textcolor{red}{(failed)}} \\ \hline

  \tabincell{c}{Exam. 78 lp\_stocfor3 \\ (n = 40216)}
  & \tabincell{c}{54 \\ (11.08)} &\tabincell{c}{8.52e-07}
  & \tabincell{c}{21+1623+50 \\ (5.69)} &\tabincell{c}{2.78e-12 }
  &\tabincell{c}{\textcolor{red}{2+4776+135} \\ \textcolor{red}{(1169.09)}} &\tabincell{c}{\textcolor{red}{5.14e+01}  \\ \textcolor{red}{(failed)}} \\ \hline

\end{tabular}}%
\end{table}%

\vskip 2mm

\vskip 2mm

\vskip 2mm

\newpage
\begin{table}[htbp]
  \newcommand{\tabincell}[2]{\begin{tabular}{@{}#1@{}}#2\end{tabular}}
  \scriptsize
  \centering
  \fontsize{8}{8}\selectfont
  \caption{Numerical results of small-scale dense LCPs (no. 1-27).}
  \label{TABDLCP1}
  \resizebox{\textwidth}{!}{
  \begin{tabular}{|c|c|c|c|c|c|c|c|c|c|}
  \hline
  \multirow{2}{*}{Problems} & \multicolumn{2}{c|}{RPFMTr} & \multicolumn{2}{c|}{PATH} & \multicolumn{2}{c|}{MILES}  \\ \cline{2-7}
              \tabincell{c}{(n*n)}
                & \tabincell{c}{steps \\(time)} & \tabincell{c}{Terr\\}
                & \tabincell{c}{major+minor+grad \\(time)} & \tabincell{c}{Terr\\}
                & \tabincell{c}{major+pivots+refactor \\(time)} & \tabincell{c}{Terr\\} \\ \hline

  \tabincell{c}{Exam. 1 lp\_25fv47 \\ (n = 2697)}
  & \tabincell{c}{56 \\ (16.83)}  &\tabincell{c}{5.85e-07}
  & \tabincell{c}{12+3457+22 \\ (200.42)} &\tabincell{c}{2.78e-09}
  &\tabincell{c}{\textcolor{red}{100+61000+600} \\ \textcolor{red}{(132.30)}} &\tabincell{c}{\textcolor{red}{3.56e+02}  \\ \textcolor{red}{(failed)}} \\ \hline

  \tabincell{c}{Exam. 2 lp\_adlittle \\ (n = 194)}
  & \tabincell{c}{43 \\ (0.26)} &\tabincell{c}{7.90e-07}
  & \tabincell{c}{8+203+10 \\ (0.03)} &\tabincell{c}{1.33e-11 }
  & \tabincell{c}{1+327+3 \\ (0.02)} &\tabincell{c}{1.38e-12} \\ \hline

  \tabincell{c}{Exam. 3 lp\_afiro  \\ (n = 78)}
  & \tabincell{c}{40 \\ (0.03)} &\tabincell{c}{7.49e-07}
  & \tabincell{c}{9+73+11 \\ (0.03)} &\tabincell{c}{2.58e-11 }
  & \tabincell{c}{1+77+2 \\ (0.08)} &\tabincell{c}{2.22e-14} \\ \hline

  \tabincell{c}{Exam. 4 lp\_agg \\ (n = 1103)}
  & \tabincell{c}{53 \\ (2.78)}  &\tabincell{c}{7.48e-07}
  & \tabincell{c}{12+1284+17 \\ (5.20)} &\tabincell{c}{5.24e-08 }
  &\tabincell{c}{\textcolor{red}{100+26000+600} \\ \textcolor{red}{(47.64)}} &\tabincell{c}{\textcolor{red}{1.43e+00}  \\ \textcolor{red}{(failed)}}  \\ \hline

  \tabincell{c}{Exam. 5 lp\_agg2 \\ (n = 1274)}
  & \tabincell{c}{54 \\ (2.80)} &\tabincell{c}{7.30e-07}
  & \tabincell{c}{10+1071+13 \\ (4.61)} &\tabincell{c}{5.44e-07 }
  & \tabincell{c}{1+132+8 \\ (2.30)} &\tabincell{c}{3.58e-09} \\ \hline

  \tabincell{c}{Exam. 6 lp\_agg3 \\ (n = 1274)}
  & \tabincell{c}{44 \\ (1.64)} &\tabincell{c}{7.03e-07}
  & \tabincell{c}{10+1053+13 \\ (4.58)} &\tabincell{c}{1.51e-07 }
  & \tabincell{c}{1+1335+8 \\ (2.25)} &\tabincell{c}{3.42e-11} \\ \hline

  \tabincell{c}{Exam. 7 lp\_bandm \\ (n = 777)}
  & \tabincell{c}{51 \\ (0.78)} &\tabincell{c}{9.17e-07}
  & \tabincell{c}{14+1075+18 \\ (2.80)} &\tabincell{c}{8.96e-11 }
  &\tabincell{c}{\textcolor{red}{100+23867+901} \\ \textcolor{red}{(37.73)}} &\tabincell{c}{\textcolor{red}{2.74e+02}  \\ \textcolor{red}{(failed)}}  \\ \hline

  \tabincell{c}{Exam. 8 lp\_beaconfd \\ (n = 468)}
  & \tabincell{c}{51 \\ (0.25)} &\tabincell{c}{8.57e-07}
  & \tabincell{c}{14+646+17 \\ (0.45)} &\tabincell{c}{9.99e-07 }
  &\tabincell{c}{\textcolor{red}{100+25300+500} \\ \textcolor{red}{(4.72)}} &\tabincell{c}{\textcolor{red}{1.76e+02}  \\ \textcolor{red}{(failed)}}  \\ \hline

  \tabincell{c}{Exam. 9 lp\_blend \\ (n = 188)}
  & \tabincell{c}{44 \\ (0.03)} &\tabincell{c}{5.39e-07}
  & \tabincell{c}{13+258+15 \\ (0.03)} &\tabincell{c}{4.73e-08 }
  & \tabincell{c}{1+547+4 \\ (0.03)} &\tabincell{c}{1.83e-09} \\ \hline

  \tabincell{c}{Exam. 10 lp\_bnl1 \\ (n = 2229)}
  & \tabincell{c}{52 \\ (9.88)} &\tabincell{c}{6.71e-07}
  & \tabincell{c}{13+3083+17 \\ (56.92)} &\tabincell{c}{1.65e-09 }
  &\tabincell{c}{\textcolor{red}{100+26600+800} \\ \textcolor{red}{(121.78)}} &\tabincell{c}{\textcolor{red}{7.79e+01}  \\ \textcolor{red}{(failed)}}  \\ \hline

  \tabincell{c}{Exam. 11 lp\_bore3d \\ (n = 567)}
  & \tabincell{c}{50 \\ (0.40)} &\tabincell{c}{9.26e-07}
  & \tabincell{c}{13+745+18 \\ (0.81)} &\tabincell{c}{8.51e-08 }
  & \tabincell{c}{2+1525+12 \\ (0.38)} &\tabincell{c}{1.05e-11} \\ \hline

  \tabincell{c}{Exam. 12 lp\_brandy \\ (n = 523)}
  & \tabincell{c}{53 \\ (0.34)} &\tabincell{c}{5.13e-07}
  & \tabincell{c}{13+697+16 \\ (0.53)} &\tabincell{c}{2.73e-07 }
  &\tabincell{c}{\textcolor{red}{100+71600+600} \\ \textcolor{red}{(11.34)}} &\tabincell{c}{\textcolor{red}{2.72e+01}  \\ \textcolor{red}{(failed)}} \\ \hline

  \tabincell{c}{Exam. 13 lp\_capri \\ (n = 753)}
  & \tabincell{c}{52 \\ (0.75)} &\tabincell{c}{5.41e-07}
  & \tabincell{c}{10+761+19 \\ (1.45)} &\tabincell{c}{2.37e-09 }
  &\tabincell{c}{\textcolor{red}{100+79200+800} \\ \textcolor{red}{(37.36)}} &\tabincell{c}{\textcolor{red}{2.10e+01}  \\ \textcolor{red}{(failed)}} \\ \hline

  \tabincell{c}{Exam. 14 lp\_czprob \\ (n = 4491)}
  & \tabincell{c}{44 \\ (42.87)} &\tabincell{c}{8.31e-07}
  & \tabincell{c}{9+3102+16 \\ (650.25)} &\tabincell{c}{3.36e-07 }
  &\tabincell{c}{\textcolor{red}{100+6200+500} \\ \textcolor{red}{(120.50)}} &\tabincell{c}{\textcolor{red}{5.33e+00}  \\ \textcolor{red}{(failed)}}  \\ \hline

  \tabincell{c}{Exam. 15 lp\_degen2 \\ (n = 1201)}
  & \tabincell{c}{49 \\ (2.32)} &\tabincell{c}{7.63e-07}
  & \tabincell{c}{12+1582+14 \\ (7.38)} &\tabincell{c}{2.83e-07 }
  & \tabincell{c}{3+566+16 \\ (2.25)} &\tabincell{c}{2.47e-09} \\ \hline

  \tabincell{c}{Exam. 16 lp\_degen3 \\ (n = 4107)}
  & \tabincell{c}{45 \\ (35.17)} &\tabincell{c}{6.85e-07}
  & \tabincell{c}{9+4003+20 \\ (388.39)} &\tabincell{c}{9.91e-07 }
  &\tabincell{c}{\textcolor{red}{100+66800+600} \\ \textcolor{red}{(129.42)}} &\tabincell{c}{\textcolor{red}{9.50e+01}  \\ \textcolor{red}{(failed)}}  \\ \hline

  \tabincell{c}{Exam. 17 lp\_e226 \\ (n = 695)}
  & \tabincell{c}{55 \\ (0.65)} &\tabincell{c}{6.15e-07}
  & \tabincell{c}{14+953+17 \\ (1.14)} &\tabincell{c}{1.45e-10 }
  &\tabincell{c}{\textcolor{red}{100+169800+1100} \\ \textcolor{red}{(34.97)}} &\tabincell{c}{\textcolor{red}{3.32e+01}  \\ \textcolor{red}{(failed)}}  \\ \hline

  \tabincell{c}{Exam. 18 lp\_etamacro \\ (n = 1216)}
  & \tabincell{c}{51 \\ (2.46)} &\tabincell{c}{8.41e-07}
  & \tabincell{c}{7+640+18 \\ (6.08)} &\tabincell{c}{5.85e-07 }
  & \tabincell{c}{2+120+9 \\ (1.17)} &\tabincell{c}{1.36e-08} \\ \hline

  \tabincell{c}{Exam. 19 lp\_fffff800 \\ (n = 1552)}
  & \tabincell{c}{70 \\ (5.76)} &\tabincell{c}{6.93e-07}
  & \tabincell{c}{14+2172+21 \\ (15.77)} &\tabincell{c}{2.09e-08 }
  &\tabincell{c}{\textcolor{red}{100+99000+700} \\ \textcolor{red}{(58.69)}} &\tabincell{c}{\textcolor{red}{1.09e+05}  \\ \textcolor{red}{(failed)}}  \\ \hline

  \tabincell{c}{Exam. 20 lp\_finnis \\ (n = 1561)}
  & \tabincell{c}{46 \\ (3.87)} &\tabincell{c}{9.59e-07}
  & \tabincell{c}{10+1583+18 \\ (14.41)} &\tabincell{c}{1.92e-07 }
  &\tabincell{c}{\textcolor{red}{100+24400+500} \\ \textcolor{red}{(28.81)}} &\tabincell{c}{\textcolor{red}{5.84e+01}  \\ \textcolor{red}{(failed)}}  \\ \hline

  \tabincell{c}{Exam. 21 lp\_fit1d \\ (n = 1073)}
  & \tabincell{c}{54 \\ (2.09)} &\tabincell{c}{8.83e-07}
  & \tabincell{c}{8+681+13 \\ (0.53)} &\tabincell{c}{3.46e-10 }
  & \tabincell{c}{1+152+4 \\ (0.06)} &\tabincell{c}{5.64e-09} \\ \hline

  \tabincell{c}{Exam. 22 lp\_ganges \\ (n = 3015)}
  & \tabincell{c}{40 \\ (14.89)} &\tabincell{c}{5.23e-07}
  & \tabincell{c}{8+1910+28 \\ (111.67)} &\tabincell{c}{6.15e-08 }
  &\tabincell{c}{\textcolor{red}{82+104468+1476} \\ \textcolor{red}{(1003.45)}} &\tabincell{c}{\textcolor{red}{9.53e+00}  \\ \textcolor{red}{(failed)}} \\ \hline

  \tabincell{c}{Exam. 23 lp\_gfrd\_pnc \\ (n = 1776)}
  & \tabincell{c}{64 \\ (7.15)} &\tabincell{c}{6.02e-07}
  & \tabincell{c}{14+1694+20 \\ (40.06)} &\tabincell{c}{2.69e-09 }
  &\tabincell{c}{\textcolor{red}{100+6800+300} \\ \textcolor{red}{(14.97)}} &\tabincell{c}{\textcolor{red}{2.16e+03}  \\ \textcolor{red}{(failed)}} \\ \hline

  \tabincell{c}{Exam. 24 lp\_grow15 \\ (n = 945)}
  & \tabincell{c}{42 \\ (1.17)} &\tabincell{c}{7.34e-07}
  & \tabincell{c}{8+856+10 \\ (2.91)} &\tabincell{c}{9.04e-13 }
  & \tabincell{c}{1+1+11 \\ (1.02)} &\tabincell{c}{1.44e-09} \\ \hline

  \tabincell{c}{Exam. 25 lp\_grow22 \\ (n = 1386)}
  & \tabincell{c}{41 \\ (2.65)} &\tabincell{c}{9.94e-07}
  & \tabincell{c}{9+1535+11 \\ (12.94)} &\tabincell{c}{4.57e-08 }
  &\tabincell{c}{\textcolor{red}{100+2600+400} \\ \textcolor{red}{(18.05)}} &\tabincell{c}{\textcolor{red}{1.68e+00}  \\ \textcolor{red}{(failed)}} \\ \hline

  \tabincell{c}{Exam. 26 lp\_lotfi \\ (n = 519)}
  & \tabincell{c}{46 \\ (0.29)} &\tabincell{c}{6.47e-07}
  & \tabincell{c}{11+408+18 \\ (0.48)} &\tabincell{c}{3.33e-11 }
  &\tabincell{c}{\textcolor{red}{100+47400+500} \\ \textcolor{red}{(13.56)}} &\tabincell{c}{\textcolor{red}{4.44e-02}  \\ \textcolor{red}{(failed)}} \\ \hline

  \tabincell{c}{Exam. 27 lp\_maros \\ (n = 2812)}
  & \tabincell{c}{69 \\ (21.44)} &\tabincell{c}{8.93e-07}
  & \tabincell{c}{15+4178+17 \\ (193.91)} &\tabincell{c}{1.04e-08 }
  &\tabincell{c}{\textcolor{red}{100+14400+300} \\ \textcolor{red}{(35.28)}} &\tabincell{c}{\textcolor{red}{2.35e+04}  \\ \textcolor{red}{(failed)}} \\ \hline

\end{tabular}}%
\end{table}%

\begin{table}[htbp]
  \newcommand{\tabincell}[2]{\begin{tabular}{@{}#1@{}}#2\end{tabular}}
  \scriptsize
  \centering
  \fontsize{8}{8}\selectfont
  \caption{Numerical results of small-scale dense LCPs (no. 28-53).}
  \label{TABDLCP2}
  \resizebox{\textwidth}{!}{
  \begin{tabular}{|c|c|c|c|c|c|c|c|c|c|}
  \hline
  \multirow{2}{*}{Problems} & \multicolumn{2}{c|}{RPFMTr} & \multicolumn{2}{c|}{PATH} & \multicolumn{2}{c|}{MILES}  \\ \cline{2-7}
              \tabincell{c}{(n*n)}
                & \tabincell{c}{steps \\(time)} & \tabincell{c}{Terr\\}
                & \tabincell{c}{major+minor+grad \\(time)} & \tabincell{c}{Terr\\}
                & \tabincell{c}{major+pivots+refactor \\(time)} & \tabincell{c}{Terr\\} \\ \hline

  \tabincell{c}{Exam. 28 lp\_modszk1 \\ (n = 2307)}
  & \tabincell{c}{46 \\ (9.37)} &\tabincell{c}{8.17e-07}
  & \tabincell{c}{13+3281+16 \\ (89.69)} &\tabincell{c}{4.94e-07 }
  &\tabincell{c}{\textcolor{red}{100+34000+500} \\ \textcolor{red}{(33.89)}} &\tabincell{c}{\textcolor{red}{9.03e+00}  \\ \textcolor{red}{(failed)}} \\ \hline

  \tabincell{c}{Exam. 29 lp\_perold \\ (n = 2131)}
  & \tabincell{c}{60 \\ (10.21)} &\tabincell{c}{7.57e-07}
  & \tabincell{c}{16+2738+24 \\ (66.41)} &\tabincell{c}{1.33e-09 }
  &\tabincell{c}{\textcolor{red}{100+4200+400} \\ \textcolor{red}{(16.06)}} &\tabincell{c}{\textcolor{red}{2.88e+04}  \\ \textcolor{red}{(failed)}} \\ \hline

  \tabincell{c}{Exam. 30 lp\_pilot\_ja \\ (n = 3207)}
  & \tabincell{c}{76 \\ (32.29)} &\tabincell{c}{9.82e-07}
  & \tabincell{c}{18+4063+29 \\ (272.11)} &\tabincell{c}{7.42e-07 }
  &\tabincell{c}{\textcolor{red}{100+71396+779} \\ \textcolor{red}{(233.45)}} &\tabincell{c}{\textcolor{red}{5.77e+06}  \\ \textcolor{red}{(failed)}} \\ \hline

  \tabincell{c}{Exam. 31 lp\_qap8 \\ (n = 2544)}
  & \tabincell{c}{42 \\ (10.69)} &\tabincell{c}{8.04e-07}
  & \tabincell{c}{12+3423+14 \\ (116.20)} &\tabincell{c}{1.81e-07 }
  &\tabincell{c}{\textcolor{red}{100+23000+400} \\ \textcolor{red}{(29.77)}} &\tabincell{c}{\textcolor{red}{1.18e+01}  \\ \textcolor{red}{(failed)}} \\ \hline

  \tabincell{c}{Exam. 32 lp\_lp\_recipe \\ (n = 295)}
  & \tabincell{c}{53 \\ (0.09)} &\tabincell{c}{7.61e-07}
  & \tabincell{c}{10+240+15 \\ (0.16)} &\tabincell{c}{6.28e-08 }
  &\tabincell{c}{\textcolor{red}{100+5144+595} \\ \textcolor{red}{(4.11)}} &\tabincell{c}{\textcolor{red}{1.91e+01}  \\ \textcolor{red}{(failed)}} \\ \hline

  \tabincell{c}{Exam. 33 lp\_sc50a \\ (n = 128)}
  & \tabincell{c}{42 \\ (0.02)} &\tabincell{c}{7.79e-07}
  & \tabincell{c}{5+92+7 \\ (0.02)} &\tabincell{c}{2.25e-08 }
  & \tabincell{c}{1+135+2 \\ (0.02)} &\tabincell{c}{7.11e-14} \\ \hline

  \tabincell{c}{Exam. 34 lp\_scagr7 \\ (n = 314)}
  & \tabincell{c}{41 \\ (0.07)} &\tabincell{c}{8.19e-07}
  & \tabincell{c}{12+391+14 \\ (0.16)} &\tabincell{c}{1.24e-08 }
  & \tabincell{c}{2+452+8 \\ (0.11)} &\tabincell{c}{1.77e-13} \\ \hline

  \tabincell{c}{Exam. 35 lp\_scagr25 \\ (n = 1142)}
  & \tabincell{c}{41 \\ (1.72)} &\tabincell{c}{8.73e-07}
  & \tabincell{c}{9+648+19 \\ (3.58)} &\tabincell{c}{1.65e-07 }
  &\tabincell{c}{\textcolor{red}{100+53500+700} \\ \textcolor{red}{(33.66)}} &\tabincell{c}{\textcolor{red}{1.09e+01}  \\ \textcolor{red}{(failed)}}  \\ \hline

  \tabincell{c}{Exam. 36 lp\_scfxm1 \\ (n = 930)}
  & \tabincell{c}{54 \\ (1.37)} &\tabincell{c}{5.12e-07}
  & \tabincell{c}{12+1260+16 \\ (2.67)} &\tabincell{c}{1.22e-07 }
  &\tabincell{c}{\textcolor{red}{100+102400+1462} \\ \textcolor{red}{(75.61)}} &\tabincell{c}{\textcolor{red}{2.23e+02}  \\ \textcolor{red}{(failed)}}  \\ \hline

  \tabincell{c}{Exam. 37 lp\_scfxm2 \\ (n = 1860)}
  & \tabincell{c}{58 \\ (7.08)} &\tabincell{c}{8.13e-07}
  & \tabincell{c}{15+2685+19 \\ (29.86)} &\tabincell{c}{3.10e-10 }
  &\tabincell{c}{\textcolor{red}{100+104383+902} \\ \textcolor{red}{(94.95)}} &\tabincell{c}{\textcolor{red}{2.35e+02}  \\ \textcolor{red}{(failed)}}  \\ \hline

  \tabincell{c}{Exam. 38 lp\_scfxm3 \\ (n = 2790)}
  & \tabincell{c}{61 \\ (18.68)} &\tabincell{c}{9.23e-07}
  & \tabincell{c}{13+3864+17 \\ (135.97)} &\tabincell{c}{2.04e-07 }
  &\tabincell{c}{\textcolor{red}{100+27200+500} \\ \textcolor{red}{(46.27)}} &\tabincell{c}{\textcolor{red}{3.89e+02}  \\ \textcolor{red}{(failed)}} \\ \hline

  \tabincell{c}{Exam. 39 lp\_scorpion \\ (n = 854)}
  & \tabincell{c}{40 \\ (0.79)} &\tabincell{c}{7.15e-07}
  & \tabincell{c}{13+1094+15 \\ (2.64)} &\tabincell{c}{2.19-10 }
  &\tabincell{c}{\textcolor{red}{100+84572+800} \\ \textcolor{red}{(31.31)}} &\tabincell{c}{\textcolor{red}{1.38e+00}  \\ \textcolor{red}{(failed)}} \\ \hline

  \tabincell{c}{Exam. 40 lp\_shell \\ (n = 2313)}
  & \tabincell{c}{50 \\ (10.06)} &\tabincell{c}{6.00e-07}
  & \tabincell{c}{13+3371+15 \\ (85.95)} &\tabincell{c}{4.39e-08 }
  &\tabincell{c}{\textcolor{red}{100+56406+497} \\ \textcolor{red}{(45.75)}} &\tabincell{c}{\textcolor{red}{1.17e+02}  \\ \textcolor{red}{(failed)}} \\ \hline

  \tabincell{c}{Exam. 41 lp\_ship04l \\ (n = 2568)}
  & \tabincell{c}{42 \\ (10.56)} &\tabincell{c}{5.64e-07}
  & \tabincell{c}{12+4243+15 \\ (96.83)} &\tabincell{c}{1.66e-09 }
  &\tabincell{c}{\textcolor{red}{100+31600+500} \\ \textcolor{red}{(61.70)}} &\tabincell{c}{\textcolor{red}{1.31e+02}  \\ \textcolor{red}{(failed)}} \\ \hline

  \tabincell{c}{Exam. 42 lp\_ship04s \\ (n = 1908)}
  & \tabincell{c}{42 \\ (5.50)} &\tabincell{c}{7.15e-07}
  & \tabincell{c}{12+3170+15 \\ (32.59)} &\tabincell{c}{4.05e-07 }
  &\tabincell{c}{\textcolor{red}{100+20000+500} \\ \textcolor{red}{(32.59)}} &\tabincell{c}{\textcolor{red}{5.69e+01}  \\ \textcolor{red}{(failed)}} \\ \hline

  \tabincell{c}{Exam. 43 lp\_ship08s \\ (n = 3245)}
  & \tabincell{c}{43 \\ (18.97)} &\tabincell{c}{9.32e-07}
  & \tabincell{c}{12+5907+15 \\ (279.95)} &\tabincell{c}{6.91e-08 }
  &\tabincell{c}{\textcolor{red}{100+70700+700} \\ \textcolor{red}{(162.72)}} &\tabincell{c}{\textcolor{red}{2.38e+02}  \\ \textcolor{red}{(failed)}} \\ \hline

  \tabincell{c}{Exam. 44 lp\_ship12s \\ (n = 4020)}
  & \tabincell{c}{42 \\ (30.77)} &\tabincell{c}{6.94e-07}
  & \tabincell{c}{9+3805+21 \\ (333.53)} &\tabincell{c}{2.36e-10 }
  &\tabincell{c}{\textcolor{red}{77+205665+1472} \\ \textcolor{red}{(1002.61)}} &\tabincell{c}{\textcolor{red}{1.97e+02}  \\ \textcolor{red}{(failed)}} \\ \hline

  \tabincell{c}{Exam. 45 lp\_sierra \\ (n = 3962)}
  & \tabincell{c}{47 \\ (33.55)} &\tabincell{c}{8.12e-07}
  & \tabincell{c}{12+4069+16 \\ (319.80)} &\tabincell{c}{7.15e-09 }
  &\tabincell{c}{\textcolor{red}{100+5000+300} \\ \textcolor{red}{(35.28)}} &\tabincell{c}{\textcolor{red}{9.99e+04}  \\ \textcolor{red}{(failed)}} \\ \hline

  \tabincell{c}{Exam. 46 lp\_standgub \\ (n = 1744)}
  & \tabincell{c}{53 \\ (5.65)} &\tabincell{c}{6.12e-07}
  & \tabincell{c}{12+2064+16 \\ (19.86)} &\tabincell{c}{3.20e-08 }
  &\tabincell{c}{\textcolor{red}{100+94960+792} \\ \textcolor{red}{(146.17)}} &\tabincell{c}{\textcolor{red}{2.44e+03}  \\ \textcolor{red}{(failed)}} \\ \hline

  \tabincell{c}{Exam. 47 lp\_tuff \\ (n = 961)}
  & \tabincell{c}{61 \\ (1.72)} &\tabincell{c}{5.46e-07}
  & \tabincell{c}{15+1404+18 \\ (3.91)} &\tabincell{c}{1.50e-10 }
  &\tabincell{c}{\textcolor{red}{100+24300+939} \\ \textcolor{red}{(62.39)}} &\tabincell{c}{\textcolor{red}{7.78e+03}  \\ \textcolor{red}{(failed)}} \\ \hline

  \tabincell{c}{Exam. 48 lp\_wood1p \\ (n = 2839)}
  & \tabincell{c}{64 \\ (20.82)} &\tabincell{c}{6.11e-07}
  & \tabincell{c}{58+3749+64 \\ (473.97)} &\tabincell{c}{1.45e-09 }
  &\tabincell{c}{\textcolor{red}{100+2700+300} \\ \textcolor{red}{(10.44)}} &\tabincell{c}{\textcolor{red}{7.27e+05}  \\ \textcolor{red}{(failed)}} \\ \hline

  \tabincell{c}{Exam. 49 lpi\_box1 \\ (n = 492)}
  & \tabincell{c}{30 \\ (0.20)} &\tabincell{c}{6.72e-07}
  & \tabincell{c}{13+657+15 \\ (0.56)} &\tabincell{c}{7.62e-08 }
  & \tabincell{c}{2+2115+24 \\ (0.66)} &\tabincell{c}{2.67e-15} \\ \hline

  \tabincell{c}{Exam. 50 lpi\_cplex2 \\ (n = 602)}
  & \tabincell{c}{43 \\ (0.45)} &\tabincell{c}{5.29e-07}
  & \tabincell{c}{12+627+14 \\ (0.95)} &\tabincell{c}{6.53e-07 }
  &\tabincell{c}{\textcolor{red}{100+46800+798} \\ \textcolor{red}{(25.13)}} &\tabincell{c}{\textcolor{red}{1.48e+02}  \\ \textcolor{red}{(failed)}} \\ \hline

  \tabincell{c}{Exam. 51 lpi\_ex72a \\ (n = 412)}
  & \tabincell{c}{34 \\ (0.12)} &\tabincell{c}{5.52e-07}
  & \tabincell{c}{13+570+15 \\ (0.41)} &\tabincell{c}{5.36e-09 }
  & \tabincell{c}{1+867+6 \\ (0.14)} &\tabincell{c}{3.13e-14} \\ \hline

  \tabincell{c}{Exam. 52 lpi\_ex73a \\ (n = 404)}
  & \tabincell{c}{33 \\ (0.12)} &\tabincell{c}{9.00e-07}
  & \tabincell{c}{13+552+15 \\ (0.38)} &\tabincell{c}{7.35e-09}
  &\tabincell{c}{\textcolor{red}{100+6739+600} \\ \textcolor{red}{(9.94)}} &\tabincell{c}{\textcolor{red}{1.50e-03}  \\ \textcolor{red}{(failed)}} \\ \hline

  \tabincell{c}{Exam. 53 lpi\_mondou2 \\ (n = 916)}
  & \tabincell{c}{38 \\ (0.94)} &\tabincell{c}{5.86e-07}
  & \tabincell{c}{11+773+13 \\ (2.36)} &\tabincell{c}{6.92e-10 }
  & \tabincell{c}{2+1151+9 \\ (0.63)} &\tabincell{c}{1.55e-14} \\ \hline

\end{tabular}}%
\end{table}%

\begin{table}[htbp]
  \newcommand{\tabincell}[2]{\begin{tabular}{@{}#1@{}}#2\end{tabular}}
  \scriptsize
  \centering
  \fontsize{8}{8}\selectfont
  \caption{Numerical results of large-scale dense LCPs (no. 54-73).}
  \label{TABDLCP3}
  \resizebox{\textwidth}{!}{
  \begin{tabular}{|c|c|c|c|c|c|c|c|c|c|}
  \hline
  \multirow{2}{*}{Problems} & \multicolumn{2}{c|}{RPFMTr} & \multicolumn{2}{c|}{PATH} & \multicolumn{2}{c|}{MILES}  \\
  \cline{2-7} \tabincell{c}{(n*n)}
                & \tabincell{c}{steps \\(time)} & \tabincell{c}{Terr\\}
                & \tabincell{c}{major+minor+grad \\(time)} & \tabincell{c}{Terr\\}
                & \tabincell{c}{major+pivots+refactor \\(time)} & \tabincell{c}{Terr\\} \\ \hline

  \tabincell{c}{Exam. 54 lp\_80bau3b \\ (n = 14323)}
  & \tabincell{c}{45 \\ (1385.62)} &\tabincell{c}{5.38e-07}
  &\tabincell{c}{\textcolor{red}{5+8283+8} \\ \textcolor{red}{(18049.72)}} &\tabincell{c}{\textcolor{red}{3.92e+01}  \\ \textcolor{red}{(failed)}}
  &\tabincell{c}{\textcolor{red}{2+8245+54} \\ \textcolor{red}{(1432.94)}} &\tabincell{c}{\textcolor{red}{1.32e+03}  \\ \textcolor{red}{(failed)}} \\ \hline

  \tabincell{c}{Exam. 55 lp\_bnl2 \\ (n = 6810)}
  & \tabincell{c}{55 \\ (153.40)} &\tabincell{c}{8.21e-07}
  & \tabincell{c}{11+8765+18 \\ (4549.92)} &\tabincell{c}{3.93e-08 }
  &\tabincell{c}{\textcolor{red}{43+27624+298} \\ \textcolor{red}{(1009.19)}} &\tabincell{c}{\textcolor{red}{9.22e+01}  \\ \textcolor{red}{(failed)}} \\ \hline

  \tabincell{c}{Exam. 56 lp\_cre\_a  \\ (n = 10764)}
  & \tabincell{c}{61 \\ (583.43)} &\tabincell{c}{6.36e-07}
  &\tabincell{c}{\textcolor{red}{4+3835+7} \\ \textcolor{red}{(18099.73)}} &\tabincell{c}{\textcolor{red}{7.38e+01}  \\ \textcolor{red}{(failed)}}
  &\tabincell{c}{\textcolor{red}{138+49799+558} \\ \textcolor{red}{(1003.36)}} &\tabincell{c}{\textcolor{red}{1.51e+03}  \\ \textcolor{red}{(failed)}} \\ \hline

  \tabincell{c}{Exam. 57 lp\_cre\_c \\ (n = 9479)}
  & \tabincell{c}{55 \\ (373.29)} &\tabincell{c}{8.33e-07}
  &\tabincell{c}{\textcolor{red}{7+8042+10} \\ \textcolor{red}{(18008.97)}} &\tabincell{c}{\textcolor{red}{2.68e+00}  \\ \textcolor{red}{(failed)}}
  &\tabincell{c}{\textcolor{red}{140+28754+599} \\ \textcolor{red}{(1004.88)}} &\tabincell{c}{\textcolor{red}{3.55e+03}  \\ \textcolor{red}{(failed)}} \\ \hline

  \tabincell{c}{Exam. 58 lp\_cycle\\ (n = 5274)}
  & \tabincell{c}{62 \\ (94.96)} &\tabincell{c}{5.76e-07}
  & \tabincell{c}{12+7861+17 \\ (2322.98)} &\tabincell{c}{2.77e-07 }
  &\tabincell{c}{\textcolor{red}{100+15300+400} \\ \textcolor{red}{(148.55)}} &\tabincell{c}{\textcolor{red}{2.07e+03}  \\ \textcolor{red}{(failed)}}  \\ \hline

  \tabincell{c}{Exam. 59 lp\_d6cube \\ (n = 6599)}
  & \tabincell{c}{54 \\ (141.97)} &\tabincell{c}{6.46e-07}
  & \tabincell{c}{10+6621+16 \\ (725.33)} &\tabincell{c}{5.15e-07 }
  &\tabincell{c}{\textcolor{red}{80+66400+800} \\ \textcolor{red}{(1009.78)}} &\tabincell{c}{\textcolor{red}{2.19e+03}  \\ \textcolor{red}{(failed)}} \\ \hline

  \tabincell{c}{Exam. 60 lp\_fit2d \\ (n = 10549)}
  & \tabincell{c}{59 \\ (550.98)} &\tabincell{c}{7.06e-07}
  & \tabincell{c}{8+5337+13 \\ (185.28)} &\tabincell{c}{2.82e-09 }
  & \tabincell{c}{2+147+4 \\ (1.11)} &\tabincell{c}{3.78e-10} \\ \hline

  \tabincell{c}{Exam. 61 lp\_greenbea \\ (n = 7990)}
  & \tabincell{c}{50 \\ (223.83)} &\tabincell{c}{8.66e-07}
  & \tabincell{c}{12+11911+16 \\ (7877.61)} &\tabincell{c}{5.44e-09 }
  &\tabincell{c}{\textcolor{red}{100+15200+400} \\ \textcolor{red}{(373.86)}} &\tabincell{c}{\textcolor{red}{1.14e+02}  \\ \textcolor{red}{(failed)}} \\ \hline

  \tabincell{c}{Exam. 62 lp\_greenbeb \\ (n = 7990)}
  & \tabincell{c}{50 \\ (225.27)} &\tabincell{c}{9.43e-07}
  & \tabincell{c}{12+12033+16 \\ (10248.31)} &\tabincell{c}{2.46e-08 }
  &\tabincell{c}{\textcolor{red}{100+17600+500} \\ \textcolor{red}{(575.45)}} &\tabincell{c}{\textcolor{red}{1.14e+02}  \\ \textcolor{red}{(failed)}} \\ \hline

  \tabincell{c}{Exam. 63 lp\_ken\_07 \\ (n = 6028)}
  & \tabincell{c}{42 \\ (96.95)} &\tabincell{c}{8.68e-07}
  & \tabincell{c}{8+4721+28 \\ (2922.75)} &\tabincell{c}{1.13e-07 }
  &\tabincell{c}{\textcolor{red}{56+6944+448} \\ \textcolor{red}{(1008.24)}} &\tabincell{c}{\textcolor{red}{5.01e+01}  \\ \textcolor{red}{(failed)}}\\ \hline

  \tabincell{c}{Exam. 64 lp\_pds\_02 \\ (n = 10669)}
  & \tabincell{c}{48 \\ (437.80)} &\tabincell{c}{7.48e-07}
  & \tabincell{c}{8+8486+26 \\ (15408.33)} &\tabincell{c}{1.82e-08 }
  &\tabincell{c}{\textcolor{red}{100+9000+300} \\ \textcolor{red}{(504.81)}} &\tabincell{c}{\textcolor{red}{1.86e+01}  \\ \textcolor{red}{(failed)}} \\ \hline

  \tabincell{c}{Exam. 65 lp\_pilot87 \\ (n = 8710)}
  & \tabincell{c}{58 \\ (307.63)} &\tabincell{c}{5.54e-07}
  & \tabincell{c}{11+8828+23 \\ (7239.50)} &\tabincell{c}{8.07e-10 }
  &\tabincell{c}{\textcolor{red}{56+48048+560} \\ \textcolor{red}{1014.17)}} &\tabincell{c}{\textcolor{red}{9.99e+02}  \\ \textcolor{red}{(failed)}} \\ \hline

  \tabincell{c}{Exam. 66 lp\_qap12 \\ (n = 12048)}
  & \tabincell{c}{48 \\ (722.82)} &\tabincell{c}{8.45e-07}
  &\tabincell{c}{\textcolor{red}{3+5789+6} \\ \textcolor{red}{(18016.27)}} &\tabincell{c}{\textcolor{red}{2.73e+01}  \\ \textcolor{red}{(failed)}}
  &\tabincell{c}{\textcolor{red}{100+23400+400} \\ \textcolor{red}{(771.77)}} &\tabincell{c}{\textcolor{red}{2.14e+01}  \\ \textcolor{red}{(failed)}} \\ \hline

  \tabincell{c}{Exam. 67 lp\_ship08l \\ (n = 5141)}
  & \tabincell{c}{47 \\ (66.37)} &\tabincell{c}{9.12e-07}
  & \tabincell{c}{12+9419+15 \\ (1537.80)} &\tabincell{c}{1.03-09}
  &\tabincell{c}{\textcolor{red}{100+80800+900} \\ \textcolor{red}{(464.31)}} &\tabincell{c}{\textcolor{red}{4.34e+02}  \\ \textcolor{red}{(failed)}} \\ \hline

  \tabincell{c}{Exam. 68 lp\_ship12l \\ (n = 6684)}
  & \tabincell{c}{47 \\ (124.64)} &\tabincell{c}{7.67e-07}
  & \tabincell{c}{12+12138+15 \\ (7256.50)} &\tabincell{c}{1.51-07}
  &\tabincell{c}{\textcolor{red}{29+11948+145} \\ \textcolor{red}{(1001.00)}} &\tabincell{c}{\textcolor{red}{7.53e+02}  \\ \textcolor{red}{(failed)}}  \\ \hline

  \tabincell{c}{Exam. 69 lp\_truss \\ (n = 9806)}
  & \tabincell{c}{45 \\ (339.85)} &\tabincell{c}{6.61e-07}
  & \tabincell{c}{10+6890+12 \\ (10715.00)} &\tabincell{c}{2.30e-08 }
  &\tabincell{c}{\textcolor{red}{23+8685+201} \\ \textcolor{red}{(1022.70)}} &\tabincell{c}{\textcolor{red}{8.80e+00}  \\ \textcolor{red}{(failed)}} \\ \hline

  \tabincell{c}{Exam. 70 lp\_woodw \\ (n = 9516)}
  & \tabincell{c}{71 \\ (494.92)} &\tabincell{c}{7.61e-07}
  & \tabincell{c}{14+7490+23 \\ (18080.28)} &\tabincell{c}{1.35e-05 }
  &\tabincell{c}{\textcolor{red}{100+15186+299} \\ \textcolor{red}{(340.14)}} &\tabincell{c}{\textcolor{red}{2.06e+02}  \\ \textcolor{red}{(failed)}} \\ \hline

  \tabincell{c}{Exam. 71 lpi\_bgindy \\ (n = 13551)}
  & \tabincell{c}{53 \\ (1297.80)} &\tabincell{c}{7.46e-07}
  &\tabincell{c}{\textcolor{red}{1+663+6} \\ \textcolor{red}{(18181.77)}} &\tabincell{c}{\textcolor{red}{3.12e+02}  \\ \textcolor{red}{(failed)}}
  &\tabincell{c}{\textcolor{red}{9+11664+81} \\ \textcolor{red}{(1025.18)}} &\tabincell{c}{\textcolor{red}{2.08e+02}  \\ \textcolor{red}{(failed)}} \\ \hline

  \tabincell{c}{Exam. 72 lpi\_gran \\ (n = 5183)}
  & \tabincell{c}{60 \\ (81.64)} &\tabincell{c}{5.32e-07}
  & \tabincell{c}{14+5542+23 \\ (1450.42)} &\tabincell{c}{1.29e-08 }
  &\tabincell{c}{\textcolor{red}{500+145000+2500} \\ \textcolor{red}{(719.47)}} &\tabincell{c}{\textcolor{red}{2.357e+03}  \\ \textcolor{red}{(failed)}} \\ \hline

  \tabincell{c}{Exam. 73 lpi\_greenbea \\ (n = 7989)}
  & \tabincell{c}{50 \\ (211.05)} &\tabincell{c}{5.37e-07}
  & \tabincell{c}{12+12005+16 \\ (12504.91)} &\tabincell{c}{7.46e-09 }
  &\tabincell{c}{\textcolor{red}{47+15546+227} \\ \textcolor{red}{(1006.56)}} &\tabincell{c}{\textcolor{red}{1.13e+02}  \\ \textcolor{red}{(failed)}}  \\ \hline

\end{tabular}}%
\end{table}%

\newpage

\vskip 2mm

\vskip 2mm

\newpage

\section{Conclusions}

\vskip 2mm

In this paper, we give the residual regularization path-following method with the
trust-region updating strategy (RPFMTr) for the linear complementarity problem.
The new residual regularization parameter improves the robustness of the path-following
method, in comparison to the traditional complementarity regularization parameter.
Meanwhile, we prove the global convergence of the new method under the standard
assumptions without the traditional assumption condition of the priority to feasibility
over complementarity. Numerical results show that RPFMTr is a robust and efficient 
solver for the linear complementarity problem, especially for the dense case. 
Furthermore, it is more robust and faster than some state-of-the-art solvers such 
as PATH \cite{DF1995,PATH,FM2000,FM2022} and MILES
\cite{Mathiesen1985,Rutherford1995,Rutherford2022} (the built-in subroutines of
the GAMS v28.2 (2019) environment \cite{GAMS}). The computational time of RPFMTr is
about 1/3 to 1/10 of that of PATH for the dense linear complementarity problem.
Therefore, RPFMTr is an alternating solver for the linear complementarity problem
and worth exploring further for the nonlinear complementarity problem.

\section*{Acknowledgments}
This work was supported in part by Grant 61876199 from National Natural Science
Foundation of China, Grant YBWL2011085 from Huawei Technologies
Co., Ltd., and Grant YJCB2011003HI from the Innovation Research Program of Huawei
Technologies Co., Ltd..
% The authors are grateful to two anonymous referees for their comments
% and suggestions which greatly improve the presentation of this paper.

\vskip 2mm

\noindent \textbf{Conflicts of interest / Competing interests:} Not applicable.

\vskip 2mm

\noindent \textbf{Availability of data and material (data transparency):} If it is requested, we will
provide the test data.

\vskip 2mm

\noindent \textbf{Code availability (software application or custom code):} If it is requested, we will
provide the code.

\end{document}